\documentclass[12pt]{amsart}
\usepackage[mathscr]{euscript}
\usepackage{amscd}
\theoremstyle{plain}
\newtheorem{theorem}{Theorem}[section]
\newtheorem{lemma}[theorem]{Lemma}
\newtheorem{proposition}[theorem]{Proposition}
\newtheorem{corollary}[theorem]{Corollary}

\theoremstyle{definition}

\theoremstyle{remark}
\newtheorem{remark}[theorem]{Remark}

\numberwithin{equation}{section}

\newcommand{\V}{\mathcal V}
\newcommand{\D}{\mathcal D}
\newcommand{\W}{\mathcal W}
\newcommand{\Co}{\mathcal C}
\newcommand{\B}{\mathcal B}
\newcommand{\N}{\mathcal N}

\renewcommand{\P}{\mathcal P}
\newcommand{\Le}{\mathcal L}

\newcommand{\C}{\mathbb C}
\newcommand{\R}{\mathbb R}
\newcommand{\bV}{\mathbb V}
\renewcommand{\Re}{\operatorname{Re}}
\renewcommand{\Im}{\operatorname{Im}}
\newcommand{\bP}{\mathbb P}
\newcommand{\zb}{\overline z}

\newcommand{\Lb}{\overline{L}}
\newcommand{\Zt}{\widetilde Z}
\newcommand{\Pt}{\widetilde{\Phi}}
\newcommand{\xt}{\tilde x}
\newcommand{\yt}{\tilde y}
\newcommand{\st}{\tilde s}
\newcommand{\pt}{\tilde p}
\newcommand{\z}{\zeta}
\newcommand{\nz}{\langle\zeta\rangle}
\renewcommand{\tt}{\tilde t}

\newcommand{\io}{\iota}
\newcommand{\G}{\Gamma}

\newcommand{\s}{\sigma}
\newcommand{\sit}{\tilde{\sigma}}
\newcommand{\e}{\epsilon}
\renewcommand{\l}{\lambda}
\renewcommand{\k}{\kappa}
\newcommand{\pa}{\partial}
\newcommand{\spa}{\operatorname{span}}
\newcommand{\Hom}{\operatorname{Hom}}

\def\crn#1#2{{\vcenter{\vbox{
        \hbox{\kern#2pt \vrule width.#2pt height#1pt
           }
          \hrule height.#2pt}}}}
\def\intprod{\mathchoice\crn54\crn54\crn{3.75}3\crn{2.5}2}
\def\into{\mathbin{\intprod}}

\begin{document}

\title[Edge of Wedge]
{Edge of the Wedge Theory in Hypo-Analytic Manifolds}
\author[M. G. Eastwood]{Michael G. Eastwood}
\address{
  Department of Pure Mathematics\\
  University of Adelaide\\
  South AUSTRALIA 5005}
\email{meastwoo@maths.adelaide.edu.au}

\author[C. R. Graham]{C. Robin Graham}
\address{
  Department of Mathematics\\
  University of Washington \\
  Box 354350\\
  Seattle, WA 98195-4350}
\email{robin@math.washington.edu}
\maketitle

\thispagestyle{empty}

\begin{center}\parbox{310pt}{\Small{\sc Abstract.}
This paper studies microlocal regularity properties of the distributions $f$
on a strongly noncharacteristic submanifold $E$ of a hypo-analytic manifold
$M$ that arise as the boundary values of solutions on wedges in $M$ with
edge $E$.  The hypo-analytic wave-front set of $f$ in the sense of
Baouendi-Chang-Treves is constrained as a consequence of the fact that $f$
extends as a solution to a wedge.}
\end{center}\bigbreak

\renewcommand{\thefootnote}{}
\footnotetext{This research was supported by the Australian Research
Council and the University of Washington.  This support and the hospitality
of the Universities of Adelaide and Washington is gratefully acknowledged.}

\section{Introduction}\label{intro}
There is an extensive literature concerning the local extension of CR
functions from submanifolds of $\C^m$, beginning with the seminal paper of
H. Lewy \cite{l}.  Much of the theory is described in the books \cite{b}
and \cite{ber}, and is concerned with extension of CR functions defined in
a neighborhood of a given point in the submanifold.  In \cite{tu}, Tumanov
has extended his minimality criterion for wedge extendability to the
situation in which the CR function is itself defined on a wedge in the
submanifold.  His criterion gives a general condition sufficient for 
extendability to a wedge in $\C^m$, but does
not give information about the direction of this wedge.  Our
work \cite{eg2} for hypersurfaces shows that there
are interesting phenomena associated to the directions of the wedges:  the
classical edge of the wedge theorem may or may not hold for CR functions on
a Levi-indefinite hypersurface, depending on directions of the wedges
involved.  In this paper we study the extension problem from wedges in a
general setting, with the intent of providing a good description of the
directions of the wedges.

The analytic wave-front set of a function (or distribution) on $\R^m$
provides a precise microlocal description of the directions of the wedges
to which the function extends holomorphically.  An analogous microlocal
theory of hypo-analytic wave-front sets was developed in \cite{bct} to
describe extension of CR functions.  Our results are formulated in
terms of this hypo-analytic wave-front set.  The natural
setting for this theory is that of a manifold $M$ with a hypo-analytic
structure, a generalization of the intrinsic structure induced on a generic
submanifold of $\C^m$.  The CR functions are replaced by the solutions of an
involutive subbundle $\V\subset \C TM$, and the holomorphic
coordinate functions by a choice of a maximally independent set of such
solutions, fixed up to a biholomorphic change.  The local theory of these
structures is given in \cite{t}; we review the relevant parts of \cite{t}
and \cite{bct} in \S 2.

Our setting is thus a wedge $\W$ in a hypo-analytic manifold $M$.  We
assume that the edge $E$ of $\W$ is a strongly noncharacteristic
submanifold of $M$.  If $M$ is a generic CR submanifold of $\C^m$, this
means exactly that $E$ is also a generic submanifold of $\C^m$, which is
contained in $M$.  We denote by $\iota_E:E\hookrightarrow M$ the inclusion.
Particular cases of interest are the maximal case in which $E$ is an open
subset of $M$, the minimal case in which $E$ is maximally real, and the
case in which $E$ is a noncharacteristic hypersurface in $M$.  For the
latter, a wedge with edge $E$ is an open set in $M$ having boundary
$E$.  If $p\in E$, the interior of the set of tangent vectors
at $p$ to curves in $\W$ defines the direction wedge $\G_p(\W)$, a linear
wedge in $T_pM$ with edge $T_pE$.  The interaction between the geometry
of the wedge $\W$ and the involutive structure $\V$ is captured by defining
a real subbundle
$\V^E\subset \V|_E$ by
$$\V^E=\{L\in \V |_E:\Re L\in TE\}$$
and transferring
the direction wedge to $\V$ by defining $$\G^{\V}(\W)=\{L\in \V^E:\Im L
\in \G(\W)\}.$$  Set also $$\G^T(\W)=\{\Re L:L\in \G^{\V}(\W)\}.$$
Then $\G_p^T(\W)$ is a cone in $T_pE$; in general $\G_p^T(\W)$ is contained
in a proper subspace of $T_pE$.

The Levi form is a crucial ingredient in the theory.  The set $T^{\circ}_p$
of characteristic covectors at $p\in M$ is defined by
$$T^{\circ}_p = \{\s\in T^*_pM:\s (L)=0 \mbox{ for all } L\in \V_p\}.$$
The Levi form $\Le_\sigma:\V_p\times\V_p\to\C$ for $\sigma\in T_p^{\circ}$
is the Hermitian sesquilinear form
defined by $$\Le_\sigma(L_1,L_2)=\frac{1}{2i}\sigma([L_1,\overline{L_2}])$$
for any sections $L_1$, $L_2$ of~$\V$.
If $L\in \V_p$, define $L\into \Le_{\s}:\V_p\rightarrow \C$ by
$$(L\into \Le_{\s})(L')=\Le_{\s}(L',L).$$
We set
$\Le_\sigma(L)=\Le_\sigma(L,L)$ and we define the second Levi form
\begin{equation}\label{Levi2}
\Le^2_\sigma(L)=\s([L,[L,\Lb]])
\end{equation}
for $L$ a section of $\V$. If $L\into \Le_{\s} = 0$, then
$\Le^2_{\s}(L)$ depends only on $L|_p$.

Now $E$ inherits a hypo-analytic structure from that on $M$, so according
to \cite{bct}, we can consider the hypo-analytic wave-front set
$WF^E(f)\subset T^*E\setminus \{0\}$ of a solution $f\in \D'(E)$ of the
involutive structure on $E$.  Our results give constraints on $WF^E(f)$
if $f$ is the boundary value of a solution on a wedge $\W\subset M$.  When
$M$ is a generic submanifold of $\C^m$, such constraints imply that $f$ may
be written as the sum of boundary values of holomorphic functions in
certain wedges, and in favorable circumstances, as the boundary value of
a holomorphic function in a single wedge or even a full neighborhood of a
point of $E$.

If $V$ is a vector space and $C\subset V$ is a cone, we define
the polar $C^{\circ}$, a closed convex cone in $V^*\setminus \{0\}$, by
$$
C^\circ=\{\xi\in V^*\setminus\{0\}:
                 \xi(v)\geq 0\mbox{ for all } v\in C\}.
$$

Our main results can be collected as follows.

\begin{theorem}\label{main}
Let $M$ be a hypo-analytic manifold, let $E\subset M$ be a strongly
noncharacteristic submanifold, and let $\W$ be a wedge in $M$ with edge
$E$.  Suppose that $f\in \D'(E)$ is the boundary value of a solution of
$\V$ on $\W$.  Then we have:
\begin{enumerate}
\item $WF^E(f)\subset (\G^T(\W))^{\circ}$.
\item If $p\in E$ and $\s\in T^{\circ}_p$ is such that there is
$L\in\G^{\V}_p(\W)$ so that one of the following conditions holds:
\begin{enumerate}
\item $\Le_{\s}(L)<0$
\item $\Le_{\s}(L)=0$ and $L\into \Le_{\s}\neq 0$
\item $\Le_{\s}(L)=0$ and $\sqrt{3}|\Im \Le^2_{\s}(L)|< \Re \Le^2_{\s}(L)$,
\end{enumerate}
then $\iota_E^*\s \notin WF^E(f)$.
\end{enumerate}
\end{theorem}
\noindent
The second condition in (c) is interpreted in the sense that it should hold
for some extension of $L$ as a section of $\V$.  However, because of (b),
we may as well assume in (c) that $L\into \Le_{\s} = 0$, in which case
$\Le^2_{\s}(L)$ is independent of the extension.

Part 1.\ in Theorem~\ref{main} is the wedge version of Theorem 3.3 of
\cite{bct}, which states that if $E$ is maximally real and $f$ is the
restriction of a solution defined
in a full neighborhood of $E$, then $WF^E(f)\subset \iota_E^*T^{\circ}$.
It follows from 1.\ that this same conclusion can be reached assuming that
$f$ is the boundary value of solutions from two opposite wedges, which
gives a weak version
of the classical edge of the wedge theorem in this setting.
It is always the case that $\iota_E^*T^{\circ}\subset
(\G^T(\W))^{\circ}$, so 1.\ can never be used to remove characteristic
covectors, i.e.\ elements of $\iota_E^*T^{\circ}$, from $WF^E(f)$.

The conditions in Part 2.\ allow one to remove characteristic covectors
from $WF^E(f)$.  Condition (a) is the hypo-analytic wedge version of the
H. Lewy extension theorem, and generalizes Theorem 6.1 of \cite{bct} for
solutions defined in a full neighborhood of $p$.  Condition (b) is an
analogue for general hypo-analytic manifolds and wedges of a result of
\cite{eg2} for hypersurfaces in $\C^m$ and for maximally real edges.
However, the corresponding result (Remark 4.3) in \cite{eg2} assumes
instead of $L\into \Le_{\s}\neq 0$ the weaker condition that $\Le_{\s}$ is
indefinite. It would be interesting to determine if (b) holds with this
weaker hypothesis in the general case.  Condition (c) is a wedge version of
a theorem of Chang \cite{chang}, who proved the analogous result for
solutions defined in a full neighborhood under the hypotheses $\Le_{\s}= 0$
and $\Le^2_{\s}(L)\neq 0$ for some $L\in \V_p$.
Our result is stronger than Chang's even in this case, as we
replace the condition  $\Le_{\s}= 0$ by the weaker hypothesis that
$\Le_{\s}(L)=0$ for a section $L$ of $\V$ for which $\Le^2_{\s}(L)\neq 0$.

All parts of Theorem 1.1 are proved using a characterization of \cite{bct}
of the complement of the wave-front set in terms of exponential decay of an
adapted version of the FBI transform on a maximally real submanifold $X$ of
$E$.  In all cases, the decay of the FBI transform is established by a
deformation of contour corresponding geometrically to the choice of a
submanifold $Y_+\subset \W$ whose boundary is $X$.  Even though our results
are more general, our arguments are simpler
than those of \cite{bct} and \cite{chang}.  We are able to achieve this 
by systematically using a
reduction introduced in \cite{chang}, which allows one to assume that the
involutive structure of $M$ is a CR structure, and by choosing $X$ and $Y_+$
geometrically {\em before} introducing special choices of coordinates.  The
changes of coordinates used in \cite{bct} and \cite{chang} in
the proofs of the full neighborhood versions of (a) and (c) use the fact
that the solution is defined in an open set and cannot be applied in the
wedge setting.

Background information and results are given in \S 2.  In \S3,
we prove Theorem~\ref{main} and show how to use this to recover results of
\cite{eg2} on extension of CR functions from wedges in hypersurfaces of
$\C^m$.  In \S 4, we discuss the involutive structure on the blow-up of a
manifold with involutive structure along a
strongly noncharacteristic submanifold; this is a construction intimately
connected with the geometry of wedges, which provides an alternative
method of studying solutions on a wedge.

\section{Hypo-analytic Structures}\label{hypo}
We begin by summarizing some of the basic notions and main results we will
need about hypo-analytic manifolds.  Some of this material is covered in
more detail in \cite{t} and \cite{bct}.

A hypo-analytic structure on a smooth manifold $M$ of dimension $n+m$ (with
$n \geq 0$ and $m\geq 1$)
consists of a choice in a neighborhood of each point of $M$ of $m$ smooth
complex valued functions $Z_1,\dots,Z_m$ with $dZ_1,\dots,dZ_m$
everywhere linearly independent, determined in overlaps up to a local
biholomorphism of $\C^m$.  The number $m$ is sometimes called the
dimension of the hypo-analytic structure and $n$ its
codimension.  A basic example is a generic CR submanifold $M$
of $\C^m$, or more generally of a $m$-dimensional complex manifold, in
which the functions $Z_j$ are taken to be the restrictions of the complex
coordinate functions.  That $M$ is a generic CR submanifold means exactly
that $dZ_1,\dots,dZ_m$ are everywhere linearly independent when restricted
to $TM$.  A function $f$ on a hypo-analytic manifold $M$
is said to be hypo-analytic if in a neighborhood of each point $p \in M$ it
is of the
form  $f=h(Z_1,\dots,Z_m)$ for a holomorphic function $h$ defined in a
neighborhood in $\C^m$ of $((Z_1(p),\dots,Z_m(p))$.  Thus for generic
CR submanifolds of $\C^m$, the hypo-analytic functions are the restrictions
to $M$ of holomorphic functions defined in a neighborhood.

The subbundle $T' \subset \C T^*M$ spanned by $dZ_1,\dots,dZ_m$ is
independent of the choice of hypo-analytic chart $(Z_1,\dots,Z_m)$ and is
called the structure bundle
of the hypo-analytic structure.  Its annihilator
$\V \equiv {T'}^{\perp}$ is a subbundle of $\C TM$ of dimension $n$, and is
involutive in the sense that the space of its sections is closed under Lie
bracket.  The bundle $\V$ is called the involutive structure
underlying the hypo-analytic structure.
We sometimes write $T'M$ and $\V M$ for $T'$ and $\V$.
In general, an involutive structure
may underlie different hypo-analytic structures. A distribution $f$ on $M$
is said to be a solution of $\V$ if $Lf=0$ for all smooth sections
$L$ of $\V$.  The involutive structure is said to be a CR structure if
$\V\cap \overline{\V}=\{0\}$.  If a hypo-analytic
structure has involutive structure $\V$ which is CR, then the map
$Z=(Z_1,\ldots,Z_m)$ defines a
local diffeomorphism from $M$ to a generic CR submanifold of $\C^m$,
so locally the hypo-analytic structures of CR type are exactly the
structures induced on generic submanifolds.

A real cotangent vector $\s \in T_p^*M$ is said to be characteristic
for the involutive structure $\V$ if $\s (L) = 0$ for all $L \in \V_p$, and
the space of characteristic covectors at $p$ is denoted $T_p^\circ$ or
$T_p^\circ M$.
The dimension $d$ of $T_p^\circ$ need not be constant as $p$ varies, so
$T^\circ$ is
not in general a vector bundle.  However, $d$ is easily seen to be
upper-semicontinuous.

A smooth submanifold $E$ of $M$ is said to be strongly noncharacteristic
if $\C T_pM=\V_p + \C T_pE$ for each $p\in E$, and maximally real if $\C
T_pM=\V_p \oplus \C T_pE$.
We will need two basic facts concerning
solutions of $\V$ near a strongly noncharacteristic submanifold $E$ of $M$.
The first (Proposition I.4.3 of \cite{t})
is that if local coordinates for $M$ are chosen near a point of $E$ such
that $E$ is defined by the vanishing of a subset of the coordinates,
then any solution of $\V$ near $E$ is a smooth function
of the variables transverse to $E$ valued in distributions in the variables
along $E$.  In particular, the restriction of any solution of $\V$ to any
strongly noncharacteristic submanifold is well-defined.  The second fact
(Corollary II.3.7 of \cite{t}) is the
uniqueness result that if a solution of $\V$ vanishes when restricted to a
strongly noncharacteristic submanifold $E$, then it must vanish in a
neighborhood of $E$.

Observe that if $E$ is a strongly noncharacteristic submanifold of a
hypo-analytic manifold $M$, there is an induced hypo-analytic structure on
$E$ since $dZ_1,\dots,dZ_m$ remain linearly independent when restricted to
$TE$.  We denote by $\V E= \V M \cap \C TE$ the induced involutive
structure, which for $p\in E$ satisfies $\dim_{\C}\V_p E = \dim E -m$.  It
is a consequence of the uniqueness result above that a (distribution)
solution  of the involutive structure on $M$ is hypo-analytic near a point
$p$ of $E$ if and only if its restriction to $E$ is
hypo-analytic near $p$ for the induced structure.
If $X\subset M$ is a maximally real submanifold and $p\in X$,
then $\V_p X=\{0\}$; whence every distribution on $X$ is a solution.

If $E\subset M$ is a strongly noncharacteristic submanifold and $p\in E$,
define
$$\V_p^E=\{L\in \V_p M: \Re L \in T_pE\}.$$
We sometimes write $\V_p^E M$ for $\V_p^E$.

\begin{lemma}\label{rank}
$\V^E$ is a real subbundle of $\V M|_E$ of
rank $=n+\dim E -m$.  The map $\Im$ which takes the imaginary part induces
an isomorphism $\V^E/\V E \cong TM|_E/TE$.
\end{lemma}

\begin{proof}
For $p\in E$, the map $\Im:\V_p^E \rightarrow T_pM$
induces a map $\V_p^E \rightarrow T_pM/T_pE$ which we claim
is surjective.  In fact, if $v\in T_pM$, then since $E$ is strongly
noncharacteristic, we can find $L\in \V_pM$ and $w\in \C T_pE$ so that
$iv=L+w$.  Taking the real part shows that $L\in \V_p^E$, whereupon taking
the imaginary part shows that $v=\Im L +\Im w$ as desired.  The kernel of
the map $\V_p^E \rightarrow T_pM/T_pE$ is $\V_pE$, so there is an induced
isomorphism $\V_p^E/\V_pE\cong T_pM/T_pE$.  Hence
$\dim \V_p^E = \dim T_pM -\dim T_pE +\dim_{\R}\V_pE=n+\dim E -m$ is
independent of $p$, from which the lemma follows.
\end{proof}

Observe that $\Re \V_p^E=(\Re \V_pM)\cap T_pE$.  In general, the dimension
of this space as well as that of $\Im \V_p^E\subset T_pM$ may vary with
$p\in E$.  However, if the involutive structure of $M$ is CR, then
$\Re$ and $\Im$ are injective on $\V M$, so also on $\V^E$.
Also, for general $M$, if $X\subset M$ is maximally real then
$\V_pX=\{0\}$, so
Lemma \ref{rank} shows that in this case, $\Im$ defines an
isomorphism from $\V_p^X$ to a $n$-dimensional subspace $\N_p$ of $T_pM$
which is a canonical complement to $T_pX$ in the sense that
\begin{equation}\label{complement}
T_pM=T_pX\oplus \N_p.
\end{equation}

We next show that there are special coordinates in a neighborhood of a point
$p$ of a maximally real submanifold $X$ in which both $X$ and the basic
solutions $Z_1,\dots,Z_m$ have particularly nice representations.  To
high order at $p$, these coordinates give a fine embedding for the
structure (in the sense of \cite{t}) adapted to $X$.

Let then $X \subset M$ be a maximally real submanifold and $p \in X$.
Dual to
(\ref{complement}) is the splitting $T_p^*M=(T_pX)^\perp\oplus
\N_p^{\perp}$. One has
$$\N_p^{\perp}=
  \{\Im\beta:\beta \in T_p'\mbox{ and } \Re\beta \in (T_pX)^\perp\}$$
and a linear map $J:\N_p^\perp\rightarrow (T_pX)^\perp$
is defined by $J(\Im \beta)=-\Re\beta$ for $\beta \in T_p'$.
Clearly, $T_p^\circ\subseteq \N_p^\perp$. Set
$d=\dim_\R T_p^\circ $ and $\nu = m-d$.  Choose a basis
$\{\sigma_1,\ldots,\sigma_d\}$ for $T_p^\circ $ and extend to a basis
$\{\sigma_1,\ldots,\sigma_d,\xi_1,\ldots,\xi_\nu\}$
for~$\N_p^\perp$. Set
$\eta_j=J(\xi_j)$. Then $\{\eta_1,\ldots,\eta_\nu\}$ is linearly
independent in $(T_pE)^\perp$ and so we may extend to a basis
$\{\eta_1,\ldots,\eta_\nu,\tau_1,\ldots,\tau_{n-\nu}\}$ of~$(T_pE)^\perp$.
Thus, we have a constructed a basis
$\{\xi_j,\eta_j,\sigma_k,\tau_l\}$ for $T_p^*M$
in which one has
$$
\begin{array}{ll}T_p^\circ=\spa_\R\{\sigma_k\}&\N_p=
\{\sigma_k=\xi_j=0\}\\[4pt]
T_p' = \spa_\C\{\sigma_k,\xi_j+i\eta_j\}&
T_pE = \{\eta_j=\tau_l=0\}.\end{array}
$$

\begin{proposition}\label{coords}
Let $X$ be a maximally real submanifold of a hypo-analytic manifold $M$
and let $p \in X$.  For each integer $N>1$ there is a neighborhood $U$ of
$p$ in $M$ with local coordinates $(x_j, y_j, s_k, t_l)$ for $1 \leq j
\leq \nu$, $1 \leq k \leq d$, $1 \leq l \leq n- \nu$, and solutions
$Z_1, \dots, Z_m$ for the hypo-analytic structure on $U$, so that
$(x,y,s,t)=0$ at $p$ and so that:
$$
X \cap U = \{y = t = 0\}
$$
\begin{equation}\label{solns1}
Z_j = x_j +iy_j + i \Psi_j (x,s) \qquad 1 \leq j \leq \nu
\end{equation}
\begin{equation}\label{solns2}
Z_{\nu + k} = s_k +i \Phi_k (x,y,s,t)  \qquad 1 \leq k \leq d
\end{equation}
where the $\Psi_j (x,s)$ and $\Phi_k (x,y,s,t)$ are smooth real functions
satisfying $d\Phi_k (0)=0$ and
\begin{equation}\label{errors}
|\Psi_j (x,s)| = O((|x| + |s|)^N), \quad
|\Phi_k (x,0,s,0)| = O((|x| + |s|)^N).
\end{equation}
\end{proposition}

\begin{proof}
Let $\{\xi_j, \eta_j, \sigma_k, \tau_l\}$ be the basis for $T_p^*M$
chosen above.  For any solutions $Z_1,\dots,Z_m$, the fiber $T_p'$ is the
span of $\{dZ_1,\dots,dZ_m\}$.  Since also $T_p' = \spa \{\xi_j + i\eta_j,
\sigma_k\}$, by replacing the $Z$'s by a linear combination we can assume
that $dZ_j(p) = \xi_j + i\eta_j$ for $1 \leq j \leq \nu$ and
$dZ_{\nu + k}(p) = \sigma_k$ for $1 \leq k \leq d$.  Define an initial set
of coordinates by setting $x_j = \Re Z_j$ and $y_j = \Im Z_j$ for
$1 \leq j \leq \nu$, $s_k = \Re Z_{\nu + k}$ for $1 \leq k \leq d$, and by
choosing functions $t_l$ for $1 \leq l \leq n- \nu$ such that $t_l=0$ on
$X$ and $dt_l(p) = \tau_l$.  Since $d Z_{\nu + k}(p) = \sigma_k$ is real,
it follows that $d Z_{\nu + k}(p) =  d \Re Z_{\nu + k}(p) = d s_k(p)$.
Recalling that $T_pX = \{ \eta_j = \tau_l = 0\}$, it follows then that in
the coordinates $(x_j, y_j, s_k, t_l)$ we have $Z_j = x_j + iy_j$,
$Z_{\nu + k} = s_k + i \Phi_k(x,y,s,t)$, and $X = \{y_j = \Psi_j(x,s),
t_l = 0\}$ for real functions $\Phi_k(x,y,s,t)$ and $\Psi_j(x,s)$
satisfying $d\Phi_k(0) = d\Psi_j(0) = 0$.

The Taylor expansions to order $N$ at $p$ of the $\Phi_k$ and $\Psi_j$
may be written $\Phi_k(x,y,s,t)=p_k(x,y,s,t) + O((|x|+|y|+|s|+|t|)^N)$
and $\Psi_j(x,s) = q_j(x,s) + O((|x|+|s|)^N)$, where the $p_k$ and $q_j$
are real polynomials of degree at most $N-1$ in their respective variables.
Observe that on $X$ we have
\begin{equation}\label{poly}
\begin{array}l
Z_j = x_j+iy_j = x_j+iq_j(x,s)+O((|x|+|s|)^N),\\[7pt]
Z_{\nu + k} = s_k + i p_k(x,q(x,s),s,0) +O((|x|+|s|)^N).
\end{array}
\end{equation}
The polynomial map $\P:\R^m \rightarrow \C^m$ given by
$$\P(x,s)=(x+iq(x,s),s+ip(x,q(x,s),s,0))$$
may be extended to $\C^m$ simply by
allowing $x$ and $s$ to be complex.  Denote also by
$\P: \C^m \rightarrow \C^m$ this extension.
Since $dp(0)=dq(0)=0$, it follows that $d\P(0)=I$, so $\P$ is
invertible near $0$ as a map from $\C^m$ to $\C^m$.  Define new
solutions
$\Zt_1,\dots,\Zt_m$ near $p$ for the hypo-analytic structure by
$(\Zt_1,\dots,\Zt_m) = \P^{-1}(Z_1,\dots,Z_m)$.  It then follows from
(\ref{poly}) that on $X$ we have
$\Zt_j = x_j + O((|x|+|s|)^N)$ for $1 \leq j \leq \nu$ and
$\Zt_{\nu + k} = s_k + O((|x|+|s|)^N)$ for $1 \leq k \leq d$.

Define new coordinates $(\xt,\yt,\st,\tt)$ on $M$ near $p$ by
$\xt_j = \Re \Zt_j$, $\yt_j = \Im \Zt_j$,
$\st_k = \Re \Zt_{\nu + k}$, and $\tt_l = t_l$.  Then in these
coordinates we have $X = \{\yt_j = \widetilde \Psi_j(\xt,\st),\,
\tt_l = 0\}$ for functions $\widetilde \Psi_j(\xt,\st)$ satisfying
$\widetilde \Psi_j(\xt,\st)=O((|\xt|+|\st|)^N)$, and also
$\Im \Zt_{\nu + k} = O((|\xt|+|\st|)^N)$ on $X$.  Finally, replace
$\yt$ by $\yt + \widetilde \Psi(\xt,\st)$ but leave
$\xt,\st,\tt$ unchanged.
In these new coordinates we have $X = \{\yt = \tt = 0\}$
and the solutions $\Zt$ have the desired form.
\end{proof}

\begin{remark}
In \cite{bct}, coordinates are used which
satisfy the conditions of Proposition \ref{coords} but also for which
$\Psi_j=0$, $1 \leq j \leq \nu$ (see II (3.9), (3.10) of
\cite{bct}).  However, such coordinates do not exist in general--for
example, if the structure is complex, the existence
of such coordinates implies that $X$ is real-analytic.  It is
possible to correct the proofs of the Theorems in \cite{bct} by using
instead the coordinates given in Proposition \ref{coords}.
\end{remark}

Let $E$ be a submanifold of a smooth manifold $M$.
In a neighborhood of a point of $E$ we may introduce coordinates
$(x',x'')$ for $M$ with $x'\in \R^r$ and $x''\in \R^s$
in which $E=\{x''=0\}$.
By a wedge in $M$ with edge $E$ we will mean an open set
$\W\subset M$
which in some such coordinate system is of the form
$\W=\B\times\Co$, where $\B$ is a ball in $\R^r$ and
$\Co\subset \R^s$ is the intersection of a ball about the origin with an
open convex cone in $\R^s\setminus \{0\}$.  Of course, this representation
of $\W$ is
only local and depends on the choice of coordinates, but we are
interested in local properties near a point of $E$ and a
direction of $\W$, so this will suffice for our purposes.
If $p\in E$, we define the direction wedge $\G_p(\W)\subset T_pM$
to be the interior of
$\{c'(0)\,|\,
c :[0,1) \rightarrow M$ is a smooth curve
in $M$ satisfying $c(t) \in \W$ for $t>0$, $c(0)=p$\}.
Then $\G_p(\W)$ is a linear wedge in $T_pM$ with edge $T_pE$,
and is determined by its image in $T_pM/T_pE$, an open convex cone.  We set
$\G(\W)=\cup_{p\in E}\G_p(\W)$.

Suppose now that $M$ has an involutive structure $\V$ and that $\W$ is a
wedge in $M$ whose edge $E$ is a strongly noncharacteristic submanifold of
$M$.  According to Lemma~\ref{rank}, $\Im$ induces an isomorphism
$:\V^E/\V E\rightarrow TM|_E/TE$.  Since the direction wedge
$\G_p(\W)$ is
determined by its image in $T_pM/T_pE$, we can use $\Im$ to define
a corresponding wedge in $\V_p^E$ which carries precisely the same
information as $\G_p(\W)$.  For $p\in E$, define
$$\G_p^{\V}(\W)=\{L\in \V_p^E: \Im L\in \G_p(\W)\}.$$
Then $\G_p^{\V}(\W)$ is a linear wedge in $\V_p^E$ with edge $\V_pE$.
We also define
$$\G_p^T(\W)=\{\Re L: L\in \G_p^{\V}(\W)\}.$$
Since $\Re$ maps $\V^E_p$ surjectively to
$(\Re \V_pM)\cap T_pE\subset T_pE$, it follows that
$\G_p^T(\W)$ is an open cone in $(\Re \V_pM)\cap T_pE$.  We set
$\G^{\V}(\W)=\cup_{p\in E}\G_p^{\V}(\W)$ and $\G^T(\W)=\cup_{p\in
E}\G_p^T(\W)$.

Let $u\in\D'(\W)$ be a solution of $\V$ and let $f\in\D'(E)$.
Let $(x',x'')$ be a coordinate system in which $E=\{x''=0\}$
and suppose that $\B$ is a ball in $\R^r$ and $\Co$ the intersection of a
ball about the origin with an
open convex cone in $\R^s\setminus \{0\}$ such that $\B\times\Co\subset
\W$.  As mentioned previously, shrinking $\B$ and $\Co$ if necessary,
$u$ defines a smooth function on $\Co$ with values in
$\D'(\B)$.  We say that
$u$ has boundary value $f$, or that $f$ is the boundary value of $u$,
and write $bu=f$,
if in each such coordinate system $(x',x'')$ and for each such $\B$ and
$\Co$, when viewed as a function on $\Co$ with values in $\D'(\B)$,
$u$ extends continuously to $\Co\cup \{0\}$ and equals $f$ at $x''=0$.
Observe that this implies that $f$ is a solution of $\V E$.
We claim that if $u$ has boundary value $f$, then
$u$ is actually $C^{\infty}$ up to $x''=0$ with values in
$\D'(\B)$.  In fact, since $E$ is strongly noncharacteristic,
one may choose sections $L_1,\ldots,L_s$ of
$\V M$ near $E$ of the form $L_j=\pa_{x''_j} +
\sum_{k=1}^r\phi_{jk}(x',x'')\pa_{x'_k}$ with $\phi_{jk}$ smooth.  Since
$L_ju=0$, it follows that
$\pa_{x''_j}u=-\sum_{k=1}^m\phi_{jk}\pa_{x'_k}u$, so each $\pa_{x''_j}u$
extends continuously up to $y=0$.  The smoothness of $u$ up to $x''=0$
follows upon iteration.

If the codimension of $E$ is 1, a wedge $\W$ with edge $E$ defines a
manifold with boundary. In Definition V.6.3 of~\cite{t}, Treves defines
a notion of distribution solution for locally integrable
structures with noncharacteristic boundary.  It follows easily from
Corollaries V.6.1, V.6.2, and V.6.3 of~\cite{t}, that if $\W$ is a wedge with
strongly noncharacteristic edge of codimension 1 in a manifold with locally
integrable involutive structure, then a
distribution solution in the sense of~\cite{t} is equivalent to a solution
whose boundary value exists in our sense.  Thus the theory in~\cite{t}
applies in our situation.  As a consequence we deduce the following.
\begin{proposition}\label{uniqueness}
If $\W$ is a
wedge with strongly noncharacteristic edge $E$ in a manifold $M$ with
locally integrable involutive structure and $u\in \D'(\W)$
is a solution of $\V$ on $\W$ with
$bu=0$ on $E$, then $u$ vanishes identically in a neighborhood
of $E$ (intersected with $\W$).
\end{proposition}\noindent
In fact, the corresponding statement for
the codimension 1 case is Corollary V.5.2 of~\cite{t} (see the last
sentence of \S V.6); the proof given there also yields an estimate on the
size of the set on which $u$ vanishes.  But the general case follows from
this, since $\W$ can be swept
out near $E$ by submanifolds of dimension $=\dim E + 1$, which inherit
locally integrable involutive structures for which $u$ restricts to be a
solution with boundary value $0$.

Next we recall the hypo-analytic wave-front set of~\cite{bct}.
To begin, let $X$ be a manifold of dimension $m$ with a hypo-analytic
structure of codimension $0$; $X$ will often arise as
a maximally real submanifold in a larger hypo-analytic manifold.  If $p\in
X$ and $Z=(Z_1,\ldots,Z_m)$ is a hypo-analytic chart in $X$ near $p$,
then $Z$ is
an embedding of a neighborhood of $p$ onto a maximally real submanifold of
$\C^m$, which of course is determined up to a local biholomorphism of
$\C^m$. For the purposes of the present discussion we may identify $X$
near $p$ with its image under $Z$ endowed with the
hypo-analytic structure induced from $\C^m$.  If
$f\in\D'(X)$ and $\s\in T^*_pX\setminus \{0\}$, then $f$ is said to be
hypo-analytic at $\s$ if there are nonempty acute open convex cones
$C_1,\ldots,C_N$ in $T_pX$, satisfying $\s(v)<0$ for all $v\in C_j$,
$1\leq j\leq N$, and wedges $\W_1,\ldots,\W_N$ in $\C^m$ with edge $X$
such that $JC_j\subset \G_p(\W_j)$,
and for each $j$ there is $u_j$ holomorphic in $\W_j$ such that
$bu_j$ exists and such that $f=bu_1 + \cdots +bu_N$.
The hypo-analytic wave-front set $WF^X(f)$ of $f$ is the complement in
$T^*X \setminus \{0\}$ of the set of points at which $f$ is
hypo-analytic; it is a closed conic subset of $T^*X \setminus \{0\}$ whose
projection to $X$ is the hypo-analytic singular support of $f$.
We set $WF_p^X(f) = T^*_pX\cap WF^X(f)$.

Two results of~\cite{bct} about wave-front sets on manifolds with
a codimension 0 hypo-analytic structure are particularly relevant
for us.  The first, Theorem 2.3 of~\cite{bct}, provides a criterion for a
distribution on $X$ to be expressible as the sum of boundary values of
holomorphic functions on specified wedges.
\begin{proposition}\label{thm2.3}
Let $C_1,\ldots,C_N$ be acute open convex cones in
$T_pX$ and let
$f\in\D'(X)$.  The following two properties are equivalent:
\begin{enumerate}
\item $WF^X_p(f)\subset \cup_{j=1}^N C_j{}^{\circ}$
\item Given for each $1\leq j \leq N$ a nonempty acute open convex cone
${\tilde C}_j$ in $T_pX$
whose closure is contained in $C_j$, there are wedges $\W_j$ in
$\C^m$ with edge $X$ such that $J{\tilde C}_j\subset \G_p(\W_j)$, and
holomorphic functions $u_j$ on $\W_j$, such that $f=bu_1 + \cdots +bu_N$.
\end{enumerate}
\end{proposition}
\noindent The special case $N=1$ is especially important as it gives
a necessary and sufficient condition for $f$ to be extendible
as a holomorphic function to a single wedge with specified direction.

Crucially important for us is Theorem 2.2 of~\cite{bct}, which
gives a criterion for microlocal hypo-analyticity in terms of
the exponential decay of a suitable FBI transform.  Let $X$ be a manifold
with a hypo-analytic structure of codimension 0 as above and let $p\in X$.
We may choose our hypo-analytic chart $Z$ such that $Z(p)=0$ and
$\Im dZ(p)=0$, in which case we may take $x_j=\Re Z_j$, $1\leq j\leq m$,
as local coordinates on $X$ near $p$.  These coordinates enable us to
identify a neighborhood of $p$ in $X$ with a neighborhood of $0$ in $\R^m$
and $T^*_pX$ with $T^*_0\R^m \cong \R^m$.  Set $\Upsilon = \Im Z$ so that
\begin{equation}\label{Z}
Z(x)=x+i\Upsilon(x)
\end{equation}
and $d\Upsilon(0)=0$.
For $z,\zeta \in \C^m$ such that $|\Im \zeta|<|\Re \zeta|$, set
$\nz=(\z_1^2+\ldots +\z_m^2)^{1/2}$ and $[z]^2=z_1^2+\ldots +z_m^2$.  If
$f$ is a compactly supported distribution in a neighborhood $U$ in $\R^m$
of $0$ and $\k >0$, define
\begin{equation}\label{FBIdef}
F^{\k}(f;z,\z)=\int_U e^{-i\z\cdot Z -\k \nz [z-Z]^2} f dZ,
\end{equation}
where $dZ=dZ_1\wedge \ldots \wedge dZ_m$ and the integral is interpreted as
a distribution pairing.  Then Theorem 2.2 combined with
Remark 2.1 of~\cite{bct} can be stated as follows.
\begin{proposition}\label{FBI}
There is an absolute constant $A>0$ so that if
$\k^* = A \sup_{x\in U, |\alpha|=2} |\pa^\alpha \Upsilon(x)|$, then
the following holds.  Let $\s \in \R^m\setminus 0$ and let $f$ be a
compactly supported distribution in $U$.  Suppose that there
is a neighborhood $V$ of the origin in $\C^m$, an open cone $\Co$ in
$\C^m\setminus 0$ containing $\s$, and constants $\e,C>0$ and $\k>\k^*$
such that
\begin{equation}\label{estimate}
|F^{\k}(f;z,\z)|\leq Ce^{-\e|\z|}
\end{equation}
for all $z\in V$ and $\z\in\Co$.  Then $\s \notin WF_p^{X}(f)$.
\end{proposition}

The wave-front set for hypo-analytic structures of positive
codimension is defined in terms of that in the codimension 0 case.
Let $E$ be a hypo-analytic manifold, let $f\in \D'(E)$ be a solution of
$\V E$, and let $p\in E$.  The hypo-analytic wave-front set
$WF^E_p(f)$ is defined as a subset
of $T_p^{\circ}E\setminus \{0\}$.  Choose a maximally real submanifold $X$
of $E$ with $p\in X$, and denote by
$\iota_X:X\rightarrow E$ the inclusion.
As discussed above, the restriction $f|_X\in \D'(X)$
is defined, and also $X$ has an induced hypo-analytic structure of
codimension 0.  Therefore we may consider
the wave-front set $WF_p^X(f|_X)\subset T^*_pX$.  Since $X$ is maximally
real, $\iota_X^*:T^{\circ}E\rightarrow T^*X$ is injective.  A covector
$\s \in T_p^{\circ}E\setminus \{0\}$ is defined to be in $WF_p^E(f)$
if $\iota_X^*\s \in WF_p^X(f|_X)$.  In \cite{bct} it is shown that this
condition is independent of the chosen maximally real submanifold $X$
containing $p$, and also that for any such $X$, one has
$WF_p^X(f|_X) \subset \iota_X^*(T_p^{\circ}E)$.  Therefore, for any
maximally real $X\subset E$ containing $p$, $\iota_X^*:WF_p^E(f)\rightarrow
WF_p^X(f|_X)$ is a bijection.

We next describe a procedure of introducing new variables which we
will use, following Chang \cite{chang},
to reduce results on general hypo-analytic
manifolds to the CR case.  Let $M$ have a hypo-analytic structure of
dimension $m$ and codimension $n$ and let $p \in M$;
recall that we write $d= \dim T^\circ_p$ and $\nu=m-d$.  On a
sufficiently small neighborhood $U$ of $p$ we may choose
coordinates $(x_1,\ldots,x_m,y_1,\ldots,y_n)$ and
hypo-analytic functions $Z_1,\ldots,Z_m$ so that at $p$ we have
$dZ_j=dx_j+idy_j$, $1\leq j\leq\nu$ and $dZ_j=dx_j$,
$\nu +1\leq j\leq m$.
Set
$M'=U \times \R^{n-\nu}$, and write $(x_{m+1},\ldots,x_{m+n-\nu})$ for the
coordinates in $\R^{n-\nu}$.  The hypo-analytic functions $Z_j$ for
$1\leq j \leq m$ pull back to $M'$ to be independent of the new variables,
and we define a hypo-analytic structure on $M'$ by augmenting these $Z_j$'s
by the functions $Z_{m+l}=x_{m+l}+iy_{\nu +l}, 1\leq l\leq n-\nu$.  At any
point $p'$ of $M'$ with first component $p$, this
structure on $M'$ has $m'=m+n-\nu$, $n'=n$, $d'=d$, and $\nu'=n=n'$.
For such $p,p'$, the map $\V_pM \ni L \to L'= L
-i\sum_{l=1}^{n-\nu}(Ly_{\nu+l}) \pa_{x_{m+l}} \in \V_{p'}M'$
is an isomorphism.  It follows that any solution of $\V M$
defines a solution of $\V M'$ which is independent of the new variables.
Any characteristic covector $\s \in T^\circ_pM$ may also be regarded
as an element of $T^\circ_{p'}M' = T^\circ_pM \oplus \{0\}$,
and it is easy to see that if $L_1$ and $L_2$ are sections of $\V M$, then
$\s([L'_1,\overline{L'_2}])=\s([L_1,\overline{L_2}])$, so that the Levi
form of $\s$ on $M'$ may be identified with that on $M$.
If $E$ is a strongly noncharacteristic submanifold of $M$ containing $p$,
then $E'=E\times \R^{n-\nu}$ is a strongly noncharacteristic
submanifold of $M'$, for which $\V_{p'}^{E'}=\{L':L\in \V_p^E\}$.
If $\W$ is a wedge in $M$ with edge $E$, then
$\W'=\W \times \R^{n-\nu}$ is a wedge in $M'$ with edge $E'$, and one has
for the direction wedges, $\G_{p'}^{\V}(\W')=\{L':L\in \G_p^{\V}(\W)\}$.
If $f\in \D'(E)$, then we may view $f \in \D'(E')$ as independent of the
new variables, and directly from the definition of the hypo-analytic
wave-front set one sees that $WF_{p'}^{E'}(f)\supset WF_p^E(f)\times
\{0\}$.

We close this section with a lemma asserting that for CR structures, a
suitable maximally real submanifold can always be chosen.

\begin{lemma}\label{chooseX}
Let $M$ be a manifold with CR involutive structure $\V$.  Let $E\subset M$
be a strongly noncharacteristic submanifold, let $p\in E$, and let $L\in
\V_p^E$.  Then there is a maximally real submanifold $X\subset E$
with $p\in X$ and $L\in \V_p^X$.
\end{lemma}

\begin{proof}
We first claim that $\spa_{\C}\V_p^E=\V_pM$.  (This does not use that $\V$
is CR.)  Note that
$\spa_{\C}\V_p^E=\V^E_p+i\V_p^E$.  Using Lemma~\ref{rank}, we have
$\dim_{\R} (\V_p^E+i\V_p^E) = 2\dim \V_p^E -\dim_{\R}\V_pE
=2n$, so it must be the case that $\V^E_p+i\V_p^E=\V_pM$ as claimed.

To prove the lemma, we may assume that $L\neq 0$.
Choose a set $\{L_1,\ldots, L_n\}\subset \V_p^E$ which is a
basis over $\C$ for $\V_pM$, and for which $L_1=L$.
Next choose $V_1,\ldots, V_k \in T_pE$ such that the cosets
$\{[V_1],\ldots,[V_k]\}$ form a basis for $T_pE/(T_pE \cap \Re \V_pM)$.
(If $\V$ is CR, it will be the case that $k=d$, but it is not necessary to
argue this separately.)
Using that $\V$ is CR and $E$ strongly noncharacteristic, it is easy to
see that $\{\Re L_1,\ldots \Re L_n, V_1,\ldots, V_k\}$ forms a basis for a
maximally real subspace of $T_pE$.  We may then take for $X$ any
submanifold of $E$ near $p$ with $T_pX$ equal to this subspace.
\end{proof}

\section{Edge of the Wedge Theory}

Throughout this section we consider a distribution $f$ on a strongly
noncharacteristic submanifold $E$ in a hypo-analytic manifold $M$, such
that $f$ is
a solution of the involutive structure $\V E$ on $E$.  We denote by
$\iota_E:E\rightarrow M$ the inclusion.  As discussed in the
previous section, $E$ inherits a hypo-analytic structure
from the hypo-analytic structure on $M$.
Our first result shows that the hypo-analytic wave-front set
of $f$ is constrained if $f$ extends to a wedge as
a solution of the involutive structure of $M$.

\begin{theorem}\label{nonchar}
Let $M$ be a hypo-analytic manifold, let $E\subset M$ be a strongly
noncharacteristic submanifold, and let $\W$ be a wedge in $M$
with edge $E$.  Suppose that $f\in
\D'(E)$ is the boundary value of a solution of $\V M$ on $\W$.  Then
$WF^E(f)\subset (\G^T(\W))^\circ$.  (Here the polar refers to the duality
between $TE$ and $T^*E$.)
\end{theorem}

\begin{proof}
Let $p\in E$ and let $\s\in T^{\circ}_pE\setminus 0$ satisfy $\s \notin
(\G_p^T(\W))^\circ$; we must show that $\s\notin WF_p^E(f)$.  Construct a new
hypo-analytic manifold $M'$ by introducing new
variables as described in \S 2.  The pullback of $\s$ satisfies
at $p'\in M'$ the analogue of the condition $\s \notin (\G_p^T(\W))^\circ$,
so it suffices to prove the result on $M'$.  Therefore we may as well assume
that our original hypo-analytic structure on $M$ has involutive structure
which is CR near $p$.

Since $\s \notin (\G_p^T(\W))^\circ$, there is $L\in \G_p^{\V}(\W)$ so that
$\s(\Re L)<0$.  By Lemma~\ref{chooseX}, we can choose a maximally real
submanifold $X\subset E$ with $p\in X$ and $L\in \V_p^X$.  Necessarily we
have $\Im L\notin T_pX$.  Choose a submanifold $Y$ of $M$ such that
$X\subset Y$ and $T_pY=\spa(T_pX,\Im L)$; in particular
$\dim(Y)=\dim(X)+1$.  Since $X$ is maximally real, $Y$ inherits near $p$
a hypo-analytic structure of codimension 1 from that on $M$.  The
involutive structure of $Y$ is CR near $p$, $\V_pY$ is spanned (over $\C$)
by $L$, $X$ is a maximally real submanifold of $Y$, and $\V_p^X Y$ is
spanned (over $\R$) by $L$. We may assume that near $p$,
$Y \setminus X$ has two connected components, the one of which lies on the
side of $X$ determined by $\Im L$ we denote by $Y_+$.  Since
$\Im L\in \G_p(\W)$,
it follows that $Y_+\subset \W$ sufficiently near
$p$.  Now $Y_+$ may be regarded as a wedge in $Y$ with edge $X$.
Considering $Y$ to be the background space, we have that
$\G_p^{\V}(Y_+)$ consists of the positive multiples of $L$.

The solution on $\W$ with boundary value $f$ restricts to $Y_+$ and defines
there a solution of the involutive structure of $Y$ having boundary
value $f|_X$ on $X$.  According to the definition of the wave-front set,
the desired conclusion $\s\notin WF_p^E(f)$ is equivalent to
$\iota_X^*\s \notin WF_p^X(f|_X)$.   Since the hypothesis $\s(\Re L)<0$
holds just as well when regarding $L\in \V_pY$ and replacing $\s$ by
$\iota_X^*\s$, it follows that we can reach this same
conclusion if we prove the theorem on $Y$.

We therefore consider the situation in which $X$ is a maximally real
submanifold in the hypo-analytic manifold $Y$ whose
structure satisfies $\nu = n = 1$ at $p \in X$, and $\W=Y_+$
is one side of $X$ in $Y$.  We choose local
coordinates $(x_1,y_1,s_1,\ldots,s_d)$ for $Y$ near $p$ as
in Proposition~\ref{coords} taking $N=2$; in our situation we have
$d=m-1$.  Set
$x=(x_1,\ldots,x_m)=(x_1,s_1,\ldots,s_d)$ and $y=y_1$ and write
$\s=(\s_1,\ldots,\s_m)$; we may assume $|\s|=1$.
Renormalizing the coordinates if necessary, we may arrange that
$\W=\{y>0\}$.  Then $L$ is a positive multiple of
$\pa_{x_1}+i\pa_y$ and the hypothesis $\s(\Re L)<0$ says $\s_1<0$.
Upon relabeling, (\ref{solns1}) and (\ref{solns2}) can be written as
\begin{equation}\label{solns}
Z_1=x_1+iy+i\Phi_1(x,y),\qquad Z_{k}=x_k+i\Phi_k(x,y), \quad 2\leq k\leq m,
\end{equation}
where $\Phi_1(x,y)=\Psi(x)$ is independent of $y$.  {From}
Proposition~\ref{coords}, we know that there is $B>0$ so that
\begin{equation}\label{B}
|\Phi(x,y)|\leq B(|x|^2 +y^2)
\end{equation}
on a fixed neighborhood of the origin.

Let $U$ be a small
open ball about the origin in $\R^m$ and $\delta_0>0$ be such
that there is a solution $u$ of $\V Y$ on $U\times (0,\delta_0)$
with boundary value $f$.  We intend to
use Proposition~\ref{FBI} to show that $\s\notin WF_0(f)$.
Introduce the
$m$-form $\omega$ on $Y_+$ given by
$\omega = e^{-i\z\cdot Z -\k \nz [z-Z]^2} u dZ$.
Since $u$ is a solution, $\omega$ is closed.
Let $\phi \in C^{\infty}_0(U)$ satisfy $\phi = 1$ for $|x|\leq r$
and $\mbox{supp }\phi \subset \{|x|\leq 2r\}$, where $r>0$ will be chosen
later.  Setting $y=0$ in (\ref{solns})
shows that $Z(x,0)$ takes the form (\ref{Z}) with
$\Upsilon(x)=\Phi(x,0)$.
Define $\kappa^*$ as in Proposition~\ref{FBI} and let
$\kappa>\kappa^*$ to be determined.
Then Stokes' Theorem gives for any $0<\delta<\delta_0$,
$$
F^{\kappa}(\phi f ;z,\z)=\int_{y=0}\phi\omega
=\int_{y=\delta}\phi\omega -\int_{U\times (0,\delta)}d\phi \wedge
\omega.
$$
We shall show that if $\kappa$ is chosen sufficiently large and
$\delta,r$ sufficiently small, then
each of the two integrals on the right hand side of the above
equation satisfies an estimate of the form (\ref{estimate}).

Set $Q=\Re \{i\z\cdot Z + \kappa \langle \z\rangle[z-Z]^2\}/|\z|$
and let $Q_0$ denote $Q$ evaluated at $z=0$ and $\z=\s$.  If we show
that there is $c>0$ so that $Q_0\geq c$ for
$(x,y)\in (\mbox{supp }\phi \times \{\delta\}) \cup (\mbox{supp }d\phi
\times [0,\delta]$), then $Q\geq c/2$ for the same $(x,y)$ and for $z$ near
$0$ and $\z$ in a conic neighborhood of $\s$, and the desired estimates on
$\int_{y=\delta}\phi\omega$ and $\int_{U\times (0,\delta)}d\phi \wedge
\omega$ then follow immediately from seminorm estimates for $u$.

Now
$$Q_0 = |\s_1|y-\s\cdot \Phi + \kappa(|x|^2-y^2-2y\Phi_1-|\Phi|^2)$$
$$\geq |\s_1|y -|\Phi| +\kappa (|x|^2 -2y^2 -2|\Phi|^2),$$
which, using (\ref{B}), is
$$\geq |\s_1|y-B(|x|^2+y^2) +\kappa (|x|^2 -2y^2 -4B^2(|x|^4+y^4))$$
$$= (|\s_1|-By-2\kappa y -4\kappa B^2 y^3)y +
(\kappa -B -4\kappa B^2|x|^2)|x|^2.$$
First choose $\kappa$ so large and $r$ so small that $\kappa \geq 4B$
and  $8Br\leq 1$, and then choose $\delta$ so
small that $B\delta +2\kappa \delta +4\kappa B^2 \delta^3 \leq
|\s_1|/2$.  Then for $(x,y)\in \mbox{supp }\phi \times [0,\delta]$
we have $Q_0\geq |\s_1|y/2 + \kappa |x|^2/2$, from which we conclude that
$Q_0 \geq c=\mbox{min}(|\s_1|\delta/2, \kappa r^2/2)$ on
$(\mbox{supp }\phi \times \{\delta\}) \cup (\mbox{supp }d\phi
\times [0,\delta])$ as desired.
\end{proof}

In considering Theorem~\ref{nonchar},
recall that $\G_p^T(\W)$ is an open cone in $(\Re \V_pM) \cap T_pE$, which
in general is a proper subspace of $T_pE$.  In fact, it is easily seen that
$(\Re \V_pM) \cap T_pE = (\iota_E^*T_p^{\circ}M)^{\perp}$, where here the
$\perp$ refers to the duality between $T_pE$ and $T_p^*E$.  Since
$\iota_E^*:T^\circ M\rightarrow T^*E$ is injective, it follows that
$(\Re \V_pM) \cap T_pE$ is all of $T_pE$ if and only if
$T_p^{\circ}M=\{0\}$, i.e.\ the involutive structure on $M$ is elliptic.
Moreover, in Theorem~\ref{nonchar}, $(\G_p^T(\W))^{\circ}$
always contains $\iota_E^*T^\circ_pM\setminus \{0\}$, and
is invariant under translation by elements of $\iota_E^*T^{\circ}_pM$.
Elements of $\iota_E^*T^\circ M\setminus \{0\}$ are referred to as
the characteristic points in $T^*E$ relative to the involutive structure of
$M$.  Thus Theorem~\ref{nonchar} can never be
used to remove characteristic points from $WF^E(f)$.  If $f$ is defined
in a full neighborhood of~$p$, then
Theorem~\ref{nonchar} shows that $WF_p^E(f)\subset \iota_E^*T_p^\circ M$,
which is Theorem 3.1 of II of~\cite{bct}.  This same conclusion can also be
obtained from Theorem~\ref{nonchar} under weaker hypotheses as follows.
\begin{corollary}\label{2sidednonchar}
Let $M$ be a hypo-analytic manifold, let $X\subset M$ be a maximally real
submanifold, let $p\in X$, and let $\W^+$ and $\W^-$ be wedges in $M$ with
edge $X$ whose directions are opposite:
$\G_p(\W^+)=-\G_p(\W^-)$.  If $f\in \D'(X)$
is the boundary value of a solution of $\V M$ on $\W^+$ and also the
boundary value of a solution of $\V M$ on $W^-$, then $WF^X_p(f)\subset
\iota_X^* T^\circ_pM$.
\end{corollary}
\noindent Corollary~\ref{2sidednonchar}
follows immediately from Theorem~\ref{nonchar} and the fact
that $(\G^T_p(\W^+))^\circ \cap (\G^T_p(\W^-))^\circ \subset
\iota_X^*T^\circ_pM$.  If the structure on $M$ is elliptic, then
Corollary~\ref{2sidednonchar} reduces to the classical
edge of the wedge theorem, so may be regarded as a weak generalization of
this theorem to the hypo-analytic setting.

Our further results give conditions under which one can remove
characteristic points from the wave-front set of the boundary value of a
solution on a wedge.  The first such result is the hypo-analytic wedge
version of the classical Lewy extension theorem.

\begin{theorem}\label{negev}
Let $M$ be a hypo-analytic manifold, let $E\subset M$ be a strongly
noncharacteristic
submanifold, and let $\W$ be a wedge in $M$ with edge $E$.
Let $p \in E$ and suppose $\s \in T^\circ_pM$ satisfies $\Le_{\s}(L)<0$
for some $L\in \G_p^{\V}(\W)$.  If $f\in \D'(E)$ is the
boundary value of a solution of $\V M$ on $\W$, then $\io_E^*\s \notin
WF^E(f)$.
\end{theorem}

\begin{proof}
First construct a hypo-analytic manifold $M'$ by
introducing new variables as described in \S 2.  If we can show
that the result holds on $M'$, we can conclude that it also holds on $M$.
Therefore we may as well assume that our original hypo-analytic structure
on $M$ satisfies $\nu = n$ at $p$.

Next use Lemma~\ref{chooseX} to choose a maximally real submanifold
$X\subset E$ such that $L\in \V^X_p$.  We then construct a submanifold
$Y\supset X$ and wedge $Y_+$ as in the proof of Theorem~\ref{nonchar}.
Restriction to $T_pY$ gives an injection
$T^\circ_pM \hookrightarrow T^\circ_pY$.
If $\s \in T^\circ_pM$ and $L_1, L_2 \in \V_pY$,
then the value of the Levi form $\Le_{\s}(L_1,L_2)$ agrees when
calculated either on
$M$ or on $Y$, since $L_1, L_2$ can be extended from $p$ to remain tangent
to $Y$ on $Y$.  Thus all the hypotheses of Theorem~\ref{negev} are
satisfied on $Y$.  Since $WF^E_p(f)$ is defined in terms of $WF^X_p(f|_X)$,
if we can prove the theorem for $Y$, it follows that
it also holds for $M$.

Therefore consider the situation in which $X$ is a maximally real
submanifold in the hypo-analytic manifold $Y$ whose
structure satisfies $\nu = n = 1$ at $p \in X$, and $\W$ is one side of
$X$ in $Y$.  Choose local
coordinates $(x_1,y_1,s_1,\ldots,s_d)$ for $Y$ near $p$ as
in Proposition~\ref{coords} taking $N=4$, and such that
$\W=\{y_1>0\}$.  Then $L$ is a positive multiple of
$\pa_{\zb}=(\partial_{x_1}+i\pa_{y_1})/2$.
As in the proof of Theorem~\ref{nonchar},
set $x=(x_1,\ldots,x_m)=(x_1,s_1,\ldots,s_d)$ and $y=y_1$.
The basic hypo-analytic functions are again given by (\ref{solns}), and we
have (\ref{B}) and
\begin{equation}\label{order4}
|\Phi(x,0)| = O(|x|^4).
\end{equation}
The characteristic covector $\s$ at
$p$ is in the span of the
$dx_k$ with $k\geq 2$; upon making a real linear transformation of these
$x_k$'s and corresponding $Z_k$'s we may assume that $\s=dx_2$.

According to (\ref{order4}), the second order terms in the Taylor
expansion of $\Phi_2$ take the form
$$
\Phi_2(x,y)=ay^2 +y\sum_{k=1}^m b_k x_k +
O(|x|^3 + y^3)
$$
with $a,b_k \in \R$.  The hypo-analytic function
$$\Zt_2=Z_2-\frac{1}{2} b_1 Z_1^2 - Z_1\sum_{k=2}^m b_k Z_{k}$$
satisfies $d\Zt_2(0)=dZ_2(0)=dx_2$ and
$\Im \Zt_2 = ay^2 + O(|x|^3 + y^3)$.  If we introduce
$\xt_2=\Re \Zt_2$ and leave the remaining $x_k$'s, $Z_k$'s, and $y$
unchanged, then equations of the form (\ref{solns}) continue to hold
in the new coordinates and
(\ref{B}) and (\ref{order4}) remain valid, possibly with a different
constant $B$.
The relations
$X=\{y=0\}$ and $\s=dx_2$ still hold in the new
coordinates and $L$ is still a positive multiple of $\pa_{\zb}$.
We may therefore assume that
in the coordinates $(x,y)$, the basic hypo-analytic
functions $(Z_1,\ldots,Z_m)$ are of the form (\ref{solns})
and in addition to (\ref{B}) and (\ref{order4}), we have
for some $C>0$,
\begin{equation}\label{C}
|\Phi_2(x,y)-ay^2|\leq C(|x|^3 + y^3)
\end{equation}
on some fixed neighborhood of the origin.
It is easily seen that the vector field
\begin{equation}\label{L}
\pa_{\zb}-i\sum_{j=1}^m(\Phi_j)_{\zb}P_j
\end{equation}
is a section of $\V$ extending $\pa_{\zb}$, where
$P_j=\sum_{k=1}^m \alpha_{kj}\pa_{x_k}$ and
$(\alpha_{jk})=(\delta_{jk}+i(\Phi_j)_{x_k})^{-1}$.
An easy calculation then shows that
$\Le_{\s}(\pa_{\zb})=a/2$, so $a<0$.  Upon
rescaling $x_1$, $y$, and $Z_1$, we may assume that $a=-1$.

Once again we intend to use Proposition~\ref{FBI} to show that $\s=dx_2
\notin WF_0(f)$.  Let $U=\{|x|\leq 3r\}$, and as in
Theorem~\ref{nonchar},
choose $\phi \in C^{\infty}_0(U)$  satisfying $\phi = 1$ for $|x|\leq r$
and $\mbox{supp }\phi \subset \{|x|\leq 2r\}$.  By (\ref{order4}), the
constant $\kappa^*$ defined in Proposition~\ref{FBI} satisfies
$\kappa^*=O(r^2)$.  Choose $r$ small enough that $\kappa^*<1/4$ and
take $\kappa = 1/4$.  Just as in the proof of Theorem~\ref{nonchar}, it
suffices to show that we can choose $r$ and $\delta$ small enough that
$Q_0=\Re(i\s \cdot Z +\kappa [Z]^2)$ satisfies $Q_0\geq c>0$
for $(x,y)\in (\mbox{supp }\phi \times \{\delta\}) \cup
(\mbox{supp }d\phi \times [0,\delta])$.  Now
$$Q_0=-\Phi_2 + \kappa(|x|^2 -y^2 -2y\Phi_1-|\Phi|^2)$$
$$\geq -\Phi_2 + \kappa(|x|^2 -2y^2 -2|\Phi|^2),$$
which, using (\ref{B}) and (\ref{C}), is
$$\geq y^2-C(|x|^3+y^3) + \kappa(|x|^2-2y^2-4B^2(|x|^4+y^4))$$
$$=(1/2-Cy-B^2y^2)y^2 + (1/4-C|x|-B^2|x|^2)|x|^2.$$
If we choose $\delta$ so small that $C\delta+B^2\delta^2<1/4$
and $r$ so small that $2Cr+4B^2r^2<1/8$, then we obtain
$Q_0\geq y^2/4 +|x|^2/8$, which yields
$Q_0\geq c=\mbox{min}(\delta^2/4,r^2/8)$ on
$(\mbox{supp }\phi \times \{\delta\}) \cup
(\mbox{supp }d\phi \times [0,\delta])$ as desired.
\end{proof}

The special case in which the solution is defined in a full neighborhood of
$p$ can be obtained by taking $E$ to be an open set in $M$.  This recovers
a result of~\cite{bct}.
\begin{corollary}\label{full}
Let $M$ be a hypo-analytic manifold, let $p\in M$, let $\s\in T^\circ_p$,
and suppose that $\Le_{\s}(L)<0$ for some $L\in\V_p$.  If $u$ is a
distribution solution near $p$, then $\s\notin WF^M(u)$.
\end{corollary}

We next turn our attention to null directions for $\Le_{\s}$.
Sometimes null directions can be perturbed to directions where the Levi
form is negative.  For instance, if there is $L_0\in \V^E_p$
with $\Le_{\s}(L_0)<0$, then given a wedge $\W$ with the property that there is
$L\in \G_p^{\V}(\W)$
satisfying $\Le_{\s}(L)=0$, one can perturb $L$ to $L'\in \G_p^{\V}(\W)$
with $\Le_{\s}(L')<0$, so Theorem~\ref{negev} applies.  However, it may
happen that $\Le_{\s}$ has a negative eigenvalue as a form on
$\V_pM$, yet $\Le_{\s}(L_1,L_2)=0$ for all $L_1,L_2\in \V^E_p$; the
hypersurface $\Im z_3 = x_1 y_2-x_2y_1$ in $\C^3$ has this property for
$E=\R^3=\{\Im z=0\}$.  For this example, every direction in $\V^E$ is null,
so Theorem~\ref{negev} does not apply.  Further effort is required to
handle such situations.

Recall the second Levi form (\ref{Levi2}).  
The following is our main lemma for null directions.

\begin{lemma}\label{nullemma}
Let $Y$ be a hypo-analytic manifold of CR type whose hypo-analytic
structure is of codimension 1.  Let $X\subset Y$ be a maximally real
submanifold and let $Y_+$ be one side of $X$ in $Y$.  Let $p\in X$ and let
$\s\in T^\circ_pY$ be such that $\Le_{\s}(L)=0$ and
$\sqrt{3}|\Im \Le^2_{\s}(L)|< \Re \Le^2_{\s}(L)$, where
$L\in \G^{\V}_p(Y_+)$.
If $f\in \D'(X)$ is the boundary value of a
solution of $\V Y$ on $Y_+$, then $\iota_X^* \s\notin WF^X(f)$.
\end{lemma}

\begin{proof}
We first remark that $\G^{\V}_p(Y_+)$ consists of a single ray, so $L$ is
determined up to a positive multiple.  Also,
the condition $\Le_{\s}(L)=0$ implies $\Le_{\s}=0$, so that the
value of $\Le_{\s}^2(L)$ is well-defined
independent of the extension of $L$ from $p$.
Choose coordinates $(x,y)=(x_1,\ldots,x_m,y)$ for $Y$ near $p$ as in the
proof of Theorem~\ref{negev}.  Thus $X=\{y=0\}$ and $Y_+=\{y>0\}$,
the hypo-analytic functions $Z_1,\ldots,Z_m$ are given by (\ref{solns})
where (\ref{B}) and (\ref{order4}) hold, and $\s=dx_2$ and
$L$ is a positive multiple of $\pa_{\zb}$.
We make the same
quadratic change of variables to obtain (\ref{C}); this time the condition
$\Le_{\s}(L)=0$ implies that $a=0$.

By (\ref{order4}), the cubic expansion of $\Phi_2$ is of the form
\begin{equation}\label{cubic}
\Phi_2(x,y)=ay^3+3bx_1y^2+3cx_1^2 y +O(|x'|(|x|^2 + y^2))
+O(|x|^4 +y^4),
\end{equation}
where $x'=(x_2,\ldots,x_m)$ and $a,b,c\in \R$.  A straightforward
calculation using the extension (\ref{L}) gives
$\Le^2_{\s}(\pa_{\zb})= 2i(\Phi_2)_{z\zb\zb}(0)
=-\frac32(a+c-ib)$.  Therefore our assumption is that
$\sqrt{3}|b|<-(a+c).$

Introduce $\Zt_2=Z_2+\mu Z_1^3$, where $\mu \in \R$ is to be
determined.
If we set $\xt_2=\Re\Zt_2$ and leave $y$ and the other $x_k$ and $Z_k$
unchanged, then all of our properties continue to hold in the new
coordinates, except that in (\ref{cubic}), $a$ is replaced by $(a-\mu)$ and
$c$ by $c+\mu$.  The inequality
$\sqrt{3}|b|<-(a+c)$ is precisely the condition that one can choose $\mu$
so that the quadratic form $(a-\mu)y^2 + 3bx_1y+3(c+\mu)x_1^2$ be negative
definite.  In fact, for $\mu=(a-c)/2$, its discriminant
is $9b^2-3(a+c)^2$.  Hence, after making this change of variables and a
subsequent rescaling of $x_1$, $y$, and $Z_1$, we may assume in
(\ref{cubic}) that
\begin{equation}\label{negdef}
ay^3+3bx_1y^2+3cx_1^2 y\leq -y^3
\end{equation}
for $y\geq 0$.

Next let $\l>0$ be small and introduce new variables \vspace{1mm}\\
$\xt_1=\l^{-1}x_1, \quad \yt=\l^{-1}y,\quad
\Zt_1=\l^{-1}Z_1$\\
$\xt_2=\l^{-3}x_2,\quad \Zt_2=\l^{-3}Z_2$\\
$\xt_k=\l^{-2}x_k,\quad \Zt_k=\l^{-2}Z_k, \quad k\geq 3.$\vspace{1mm}\\
Observe first that if $(\xt,\yt)$ is in a fixed neighborhood of the origin
and $\l_0>0$ is chosen sufficiently small, then for all $\l$ with
$0<\l\leq\l_0$, the corresponding $(x,y)$
will lie in the neighborhood of the origin in which we have
been working.  This change of variables has
the effect in (\ref{solns}) of replacing $\Phi$ by new functions $\Pt$
given by \vspace{1mm}\\
$\Pt_1(\xt,\yt)=\l^{-1}\Phi_1(x,y)$\\
$\Pt_2(\xt,\yt)=\l^{-3}\Phi_2(x,y)$\\
$\Pt_k(\xt,\yt)=\l^{-2}\Phi_k(x,y),\quad k\geq 3.$\vspace{1mm}\\
We deduce that $\Pt(\xt,\yt)$ satisfies (\ref{B})
and (\ref{order4}) in the new coordinates with constants independent of
$\l$ for $0<\l\leq\l_0$.  Also, (\ref{cubic}) is replaced by an analogous
equation for $\Pt_2(\xt,\yt)$ in which the error terms
$O(|\xt'|(|\xt|^2 + \yt^2))$ and $O(|\xt|^4 +\yt^4)$
can both be bounded by $C\l (|\xt|^3 +\yt^3)$ for a
constant $C$ independent of $\l$.  Combining this observation with
(\ref{negdef}) and removing the tildes from the new coordinates, we see
that in addition to (\ref{B}) and (\ref{order4}) we can assume that
\begin{equation}\label{cubicineq}
\Phi_2(x,y)\leq -y^3 +C\l (|x|^3+y^3).
\end{equation}

In order to apply Proposition~\ref{FBI}, let $\delta>0$ to be chosen small,
let $U=\{|x|<6\delta\}$, and choose $\phi\in C^{\infty}_0(U)$ satisfying
$\mbox{supp } \phi\subset \{|x|\leq 5\delta\}$ and $\phi = 1$ if
$|x|\leq 4\delta$.  By (\ref{order4}), we have $\kappa^*=O(\delta^2)$.
Choose $\delta$ small enough to ensure that $\kappa^*<\delta/16$ and set
$\kappa = \delta/16$.  Choose $\l$ and $\delta$ small enough that
for all $(x,y)\in U\times [0,\delta]$ we have
$\l C(|x|^3+y^3)\leq \delta^3/4$ and $4B^2(|x|^4+y^4)\leq 2\delta^2$.
We shall
show that if $\delta$ is chosen small enough, then
$Q_0=\Re(i\s \cdot Z +\kappa [Z]^2)$ satisfies $Q_0\geq c>0$
for $(x,y)\in (\mbox{supp }\phi \times \{\delta\}) \cup
(\mbox{supp }d\phi \times [0,\delta])$.

As in the proof of Theorem~\ref{negev}, we have
$Q_0 \geq -\Phi_2 + \kappa(|x|^2 -2y^2 -2|\Phi|^2),$
so by (\ref{B}) and (\ref{cubicineq}) we obtain
$$Q_0 \geq y^3 -C\l(|x|^3+y^3) + \kappa(|x|^2 -2y^2 -4B^2(|x|^4+y^4))$$
$$\geq y^3 -\delta^3/4 +\kappa(|x|^2 -2y^2-2\delta^2).$$
Therefore, if $y=\delta$ and $x\in U$, then
$$Q_0\geq \delta^3 -\delta^3/4 +\kappa(-2\delta^2-2\delta^2)=\delta^3/2.$$
If $y\in [0,\delta]$ and $x\in \mbox{supp } d\phi$, then
$$Q_0\geq -\delta^3/4 +\kappa(16\delta^2-2\delta^2-2\delta^2)=\delta^3/2.$$
Thus the desired inequality holds with $c=\delta^3/2$.
\end{proof}

We present two conditions which allow one to use null
directions to remove characteristic points from the wave-front set.
The first is that the null direction is not degenerate for
the Levi form, and applies to the hypersurface $\Im z_3 = x_1 y_2-x_2y_1$
in $\C^3$ mentioned above.  For this example, every direction
$L\in\V_p^E\setminus \{0\}$ is null and nondegenerate,
so it suffices to have an extension to any wedge.

\begin{theorem}\label{nondegnull}
Let $M$ be a hypo-analytic manifold, let $E\subset M$ be a strongly
noncharacteristic submanifold, and let $\W$ be a wedge in $M$ with edge
$E$.  Let $p \in E$ and suppose $\s \in T^\circ_pM$ has the property that
for some $L \in \G_p^{\V}(\W)$, $\Le_{\s}(L)=0$ and
$L\into\Le_{\s}\neq 0$.  If $f\in \D'(E)$ is the
boundary value of a solution of $\V M$ on $\W$, then $\io_E^*\s \notin
WF^E(f)$.
\end{theorem}

\begin{proof}
By introducing extra variables we may assume that $M$ is of CR type near
$p$.

By Lemma~\ref{chooseX}, we can choose a maximally real submanifold
$X\subset E$ so that $L\in \V_p^X$.  According to the definition of the
wave-front set, it suffices to show that $\iota_X^*\s \notin WF^X(f|_X)$.
The set $\W$ is also a wedge with edge $X$, and the hypotheses of the
Theorem continue to hold if $X$ is viewed as the edge.  Therefore we may as
well prove the theorem with $E$ replaced by $X$.

If there is ${\tilde L}\in \V_p^X$ such that
$\Re \Le_{\s}(L,{\tilde L})\neq 0$, then
we can choose $\epsilon$ small and of the appropriate sign so that the
vector $L' = L+\epsilon {\tilde L}$ will satisfy $\Le_{\s}(L')<0$.
But for sufficiently small $\epsilon$, we have $L'\in \G_p^{\V}(\W)$.
Therefore in this case the conclusion follows from
Theorem~\ref{negev}. Hence we may assume that
$\Re \Le_{\s}(L,{\tilde L})= 0$ for all ${\tilde L}\in \V_p^X$.

Choose coordinates as in Proposition~\ref{coords} with $N$ large.  Since
$M$ is of CR type, we have $\nu=n$ so there are no $t$ variables.
Now $\V_p^X=\spa_{\R} \{\pa_{\zb_j}, 1\leq j \leq n\}$, so we may make a
real linear change of coordinates and corresponding change of $Z's$ to
arrange that $L=\pa_{\zb_1}$.  Similarly we may take $\s=ds_1$.
The hypothesis $L\into\Le_{\s} \neq 0$ is equivalent to
$\Le_{\s}(\pa_{\zb_1},\pa_{\zb_j})\neq 0$ for some $j\geq 2$; by reordering
the coordinates we may assume that $j=2$.
Finally, recalling our assumption that
$\Re \Le_{\s}(L,{\tilde L)}= 0$ for all ${\tilde L}\in \V_p^X$,
we conclude that
$\Le_{\s}(\pa_{\zb_1},\pa_{\zb_2})= i\lambda$ for some nonzero $\lambda \in
\R$.

It is straightforward to check using the form of the $Z's$ in
Proposition~\ref{coords} that one may choose a basis
$\{L_j, 1\leq j \leq n\}$ for $\V M$ near $0$ of the form
$$L_j = \pa_{\zb_j} + A_{jl} \pa_{s_l} + B_{jk} \pa_{x_k};$$
in fact, the equations $L_j Z_k=0$ uniquely determine the coefficients
$A_{jl}$ and $B_{jk}$, and one has that the $A_{jl}$ vanish at $0$ and
$B_{jk}=O((|z|+|s|)^{N-1})$.

Define a submanifold $Y\subset M$ by
$Y=\{(x,y,s):y_2=\alpha x_1y_1,y_3=\ldots=y_n=0\}$, where
$\alpha \in \R$ is to be determined.  Define also
$Y_+=\{(x,y,s)\in Y:y_1>0\}$.
We have that $X\subset Y$, $T_pY = \spa (T_pX, \Im L)$,
and since $\Im L\in \G_p(\W)$, it
follows that $Y_+\subset \W$ sufficiently near the origin.
$Y$ inherits a hypo-analytic structure of codimension 1 and CR type
in which $X$ is a maximally real submanifold.  At each point of $Y$,
the space $\V Y$ is
spanned by a single complex vector of the form
$L=L_1+\sum_{j=2}^n c_j L_j$; this normalization ensures that this extended
$L$ agrees with the $L$ we already have at the origin.
The $c_j$ are determined by the requirement
that $L$ should be tangent to $Y$.  The condition
$Ly_j=0$ for $3\leq j\leq n$  forces $c_j=0$ for these $j$.
This leaves the one condition
$L(y_2-\alpha x_1y_1)=0$, which is satisfied by
$L=L_1 +\phi L_2$ for $\phi=\alpha \zb_1 + O((|z|+|s|)^N)$.

A direct calculation using $\s \in T^{\circ}_p$,
$\s([L_1,\Lb_1]) = 0$, and $\s([L_1,\Lb_2]) = -2\lambda$
results in
$\Le^2_{\s}(L) = \Le^2_{\s}(L_1) +2\lambda \alpha$.
Now $\Le^2_{\s}(L_1)$ is arbitrary, but
by suitable choice of $\alpha$ we can ensure that
$\sqrt 3 |\Im \Le^2_{\s}(L)| <\Re \Le^2_{\s}(L)$.  Therefore the
result follows from Lemma~\ref{nullemma}.
\end{proof}

Now Lemma~\ref{nullemma} can be extended to general structures.

\begin{theorem}\label{chang}
Let $M$ be a hypo-analytic manifold, let $E\subset M$ be a strongly
noncharacteristic submanifold, and let $\W$ be a wedge in $M$ with edge
$E$. Let $p \in E$ and suppose $\s \in T^\circ_pM$ has the property that
there is a section $L$ of $\V M$ satisfying $L|_p \in \G_p^{\V}(\W)$,
$\Le_{\s}(L)=0$, and
$$\sqrt{3}|\Im \Le^2_{\s}(L)|< \Re \Le^2_{\s}(L).$$
If $f\in \D'(E)$ is the
boundary value of a solution of $\V M$ on $\W$, then $\io_E^*\s \notin
WF^E(f)$.
\end{theorem}

\begin{proof}
By introducing extra variables, we may assume that $M$ is of CR type near
$p$.

If $L\into\Le_{\s}\neq 0$, then the conclusion follows
from Theorem~\ref{nondegnull}.  Therefore we may assume that
$L\into\Le_{\s}=0$.   In this case, the value of $\Le^2_{\s}(L)$ is
independent of the extension of $L|_p$, so we may as well
just consider $L$ as an element of $\G_p^{\V}(\W)$.

As in the proof of Theorem~\ref{negev},
choose a maximally real submanifold $X\subset E$ such that $L\in \V^X_p$
and a submanifold $Y\subset M$ such that $X\subset Y$ and
$T_pY=\spa(T_pX,\Im L)$.
Then $L$ spans $\V^X_pY$ and satisfies $\Le_{\s}(L)=0$ and
$\sqrt{3}|\Im \Le^2_{\s}(L)|< \Re \Le^2_{\s}(L)$.
The conclusion thus follows from Lemma~\ref{nullemma}.
\end{proof}

Taking $E$ to be an open set in $M$, we obtain a strengthened version
of Chang's Theorem for solutions defined in a full neighborhood.

\begin{corollary}\label{changfull}
Let $M$ be a hypo-analytic manifold, let $p\in M$, let $\s\in T^\circ_p$
and suppose that there is a section $L$ of $\V$ near $p$ such that
$\Le_{\s}(L)=0$ and $\Le^2_{\s}(L)\neq 0$.   If $u$ is a distribution
solution near $p$, then $\s\notin WF^M(u)$.
\end{corollary}

Of course, for a given $f\in \D'(E)$, one typically uses different of the
criteria established above
at different points $\s\in T^*E\setminus \{0\}$ to constrain $WF^E(f)$.
If $M$ is a generic submanifold
of $\C^m$, one is interested in knowing whether
a CR distribution $u$ on $\W$ extends as a holomorphic function to a
particular wedge in $\C^m$.  Such results can be obtained
via Propositions~\ref{uniqueness} and \ref{thm2.3}
by constraining the wave-front set of $f=bu$ appropriately.  We illustrate
this when $M$ is a hypersurface and the edge is maximally real.
In this case $\dim T^\circ = 1$; we
denote by $\Le$ the Levi form of $M$, which is determined up to a nonzero
real multiple.

\begin{proposition}\label{onesidedhyp}
Let $M$ be a hypersurface in $\C^m$ and let
$\W\subset M$ be a wedge with maximally real edge $X$ and let $p\in X$.
Suppose that there is $L\in \G_p^{\V}(\W)$ such that
$\Le(L) =0$ and $L\into\Le\neq 0$ (in
particular, $\Le$ is indefinite).
If $u$ is a CR function on $\W$ with boundary value
$f\in \D'(X)$, then there is a wedge $\W'$ in $\C^m$ with edge $X$ such
that $L\in \G_p^{\V}(\W')$ and such that there is a
holomorphic function $u'$ on $\W'$ with $bu'=f$ and $u'=u$ on $\W'\cap\W$.
\end{proposition}

\begin{proof}
It follows from Theorem~\ref{nonchar} that $WF_p^X(f)\subset
(\G_p^T(\W))^\circ$.
Theorem~\ref{nondegnull} implies that $WF^X(f)\cap T^\circ_p M=\emptyset$.
{From} these two facts one sees easily that
there is an open acute convex cone $C\subset T_pX\setminus \{0\}$
such that $\Re L\in C$ and $WF_p^X(f)\subset C^\circ$.  By
Proposition~\ref{thm2.3}, $f$ is the boundary value of a holomorphic
function on a wedge $\W'$ in $\C^m$ with edge $X$ with $L\in \G_p^{\V}(\W')$.
By Proposition~\ref{uniqueness}, this extension agrees
with $u$ on $\W'\cap\W$.
\end {proof}

{From} Proposition~\ref{onesidedhyp}
and the classical edge of the wedge theorem for holomorphic
functions (or directly using Corollary~\ref{2sidednonchar} and
Theorem~\ref{nondegnull}), one deduces the
following edge of the wedge theorem for CR functions.

\begin{proposition}\label{twosidedhyp}
Let $M$ be a hypersurface in $\C^m$ and let
$\W^{\pm}\subset M$ be wedges with maximally real edge $X$ whose directions
at $p\in E$ are opposite:
$\G_p(\W^+)=-\G_p(\W^-)$.
Suppose that there is $L\in \G^{\V}_p(\W^+)$ such that $\Le(L) =0$ and
$L\into\Le\neq 0$.
If $u^{\pm}$ are CR functions on $\W^{\pm}$ such that $bu^+=bu^-$, then
there is a holomorphic function in a neighborhood in $\C^m$ of $p$ which
extends $u^{\pm}$.
\end{proposition}
In particular, Proposition~\ref{twosidedhyp} implies that there is a CR
function in a neighborhood in $M$ of $p$ which extends $u^{\pm}$.  Thus
this gives an intrinsic version of the classical edge of the wedge theorem
for CR functions.

Propositions \ref{onesidedhyp} and \ref{twosidedhyp} were proved in
\cite{eg2} for continuous boundary values using folding screens, however
with the condition
$L\into\Le\neq 0$ replaced by the weaker hypothesis that
$\Le$ is indefinite.  We remark that it is also possible to use
Theorem~\ref{nonchar} together with Theorem~\ref{negev} to prove a version
of Proposition \ref{onesidedhyp} if instead of a null direction there is
$L\in\G_p^{\V}(\W)$
such that $\Le(L) \neq 0$.  However, in this situation
Theorem~\ref{negev} only applies on one side of $T^\circ_pM$, so the wedge
$\W'$ obtained in $\C^m$ must be tilted slightly and need not have
$\Im L\in \G_p(\W')$.
There is also a two-sided version when $\Le(L)\neq 0$;
in that case $WF_p^X(f)$ reduces to a single characteristic ray so that the
wedge $\W'$ given by Proposition~\ref{thm2.3} is arbitrarily close (in the
conic sense) to a half space lying on one side of $T_pM$ intersected with a
neighborhood of $p$.

\section{Blow-ups}

The geometry of wedges in a manifold $M$ with a given edge $E$ is naturally
encoded in a new manifold, the blow-up $B$ of $M$ along $E$.  It turns out
that if $M$ has an involutive or hypo-analytic structure and $E$ is
strongly noncharacteristic, then $B$ inherits a structure of the same
type which can be used to study solutions on wedges in $M$.  As a natural
geometric way of constructing new structures from old, the construction may
be of more general interest.  For example, in \cite{eg1} the blow-up in the case
$\R^n\subset\C^n$ is used to prove the classical
edge of the wedge theorem, but it is also pointed out
there how the construction arises in real integral geometry. 
In this section we describe the geometry of the blow-up in the general case
and indicate its relation to the wedge regularity theorems of the previous
section.  

Let $E$ be a submanifold of a manifold $M$ and let $\Sigma\to E$ denote the
projective normal bundle of $E$ in~$M$.  As a set, the blow-up $B$ of $M$
along~$E$ is obtained by replacing $E$ with~$\Sigma$.  It is naturally a
smooth manifold containing $\Sigma$ as a hypersurface, with a smooth
blow-down mapping $b:B\to M$, 
which is the identity outside $E$ and the projection $\Sigma\to E$
over~$E$.  If local coordinates $(x,y)$, $x\in \R^k$, $y\in \R^l$, are
chosen on $M$ such that $E$ is given by $\{y=0\}$, 
then the fibers of $\Sigma$ can be identified with $\R\bP^{l-1}$, for
which we can use usual affine coordinates on a cell. The smooth structure
on $B$ near the subset of $\Sigma$ corresponding to the first cell is
defined by coordinates $(x,y_1,u)$ with $u\in \R^{l-1}$,
in which $\Sigma = \{y_1=0\}$ and the blow-down map takes the form
$$
(x,y_1,u)\stackrel {b}{\mapsto} (x,y_1,y_1 u).
$$
The image under $b$ of an open set $U$ in $B$ intersecting $\Sigma$ is a
pair of opposite wedges in
$M$ with edge $E$ and conversely.  A (one-sided) wedge in $M$ corresponds
to the intersection of $U$ with one side of $\Sigma$ in $B$.  The intersection
$U\cap \Sigma$ is the set of directions of $b(U)$.

Suppose now that $M$ has an involutive structure $\V M$.  Since
$b:B\setminus \Sigma \rightarrow M\setminus E$ is a diffeomorphism, the
involutive structure of $M$ pulls back to an involutive structure
on $B\setminus \Sigma$.  We shall show that if $E\subset M$ is strongly
noncharacteristic, then this involutive structure extends smoothly across
$\Sigma$.  According to Lemma \ref{rank}, if $p\in E$, 
then $\Im :\V^E_p/\V_pE \rightarrow T_pM/T_pE$ is an isomorphism.  
We shall use systematically this isomorphism to represent the fiber
$\Sigma_p$ as $\bP(\V^E_p/\V_pE)$.  We therefore represent a point $\pt \in
\Sigma$ as a pair $\pt = (p,\ell)$ with $p=b(\pt)$ and 
$\ell \in \bP(\V^E_p/\V_pE)$.  We
choose $L\in \V_p^E$ representing $\ell$ and write $\ell =[L]$.  Of course,
$L$ is unique only up to a nonzero real scale factor and up to addition of
an element of $\V_pE$.  

\begin{theorem}\label{involB}
Let $M$ have an involutive structure $\V M$ and let $E\subset M$ be a
strongly noncharacteristic submanifold.  Let $B$ be the blow-up of $M$
along $E$, let $\pt \in B$, and set $p=b(\pt)$.  Define
$$
\V_{\pt}B =
\left\{
\begin{array}{ll}
(b_*)^{-1}(\V_pM) &  \mbox{if} \,\,\pt\notin\Sigma \\
(b_*)^{-1}(\C L \oplus \V_pE) & \mbox{if} \,\,
\pt = (p,[L])\in\Sigma,
\end{array}
\right.
$$
where $b_*:\C T_{\pt}B\rightarrow \C T_p M$ is the differential of $b$ at 
$\pt$.  Then $\V B$ is a smooth involutive structure on $B$.
\end{theorem}

\begin{proof}
It suffices to reason near $\pt \in \Sigma$.  Elementary linear algebra 
(cf.\ the discussion prior to Proposition~\ref{coords}) shows that one can 
choose coordinates on $M$ near $p\in E$ of the form $(x, y, s, t)$ for 
$x,y\in \R^{\nu}$, $s\in \R^d$, $t\in \R^{n- \nu}$, where $d+\nu=m$,
so that $(x,y,s,t)=0$ at $p$ and so that:
$$
E = \{y_{r+1}=\ldots = y_{\nu} = t_{\rho +1}=\ldots = t_{n-\nu} = 0\} 
$$
for some $r$ and $\rho$ satisfying $0\leq r \leq \nu$,
$0\leq \rho \leq n-\nu$, and such that
\begin{equation}\label{basis}
T'_pM = \spa\{dx_1+idy_1,\ldots dx_{\nu}+idy_{\nu}, ds_1,\ldots
ds_d\}.  
\end{equation}
Write $y=(y',y'')$ with 
$y'=(y_1,\ldots ,y_r)$, $y''=(y_{r+1},\ldots ,y_{\nu})$,
and $t=(t',t'')$ with
$t'=(t_1,\ldots ,t_{\rho})$, $t''=(t_{\rho +1},\ldots ,t_{n-\nu})$. 
The corresponding coordinates on $B$ are obtained as described above.  We 
consider the case in which
$\pt$ corresponds to a line in the normal space
$\{(y'',t'')\}$ which satisfies 
$y_{r+1}\neq 0$; the argument in other cases is similar.  
Thus we introduce coordinates 
$u=(u_{r+2},\ldots ,u_{\nu})\in \R^{\nu-r-1}$,  
$v=(v_{\rho +1},\ldots,v_{n-\nu})\in \R^{n-\nu-\rho}$ 
and obtain coordinates on $B$ near $\pt$ and a formula for $b$:
\begin{equation}\label{b}
(x, y',y_{r+1},u,s,t',v) \stackrel{b}{\mapsto}
(x,y',y_{r+1},y_{r+1}u,s,t',y_{r+1}v).
\end{equation}
The point $(x, y',0,u,s,t',v)\in \Sigma$ corresponds to the direction
$[L]$ with 
\begin{equation}\label{Lform}
L=2\pa_{\zb_{r+1}}+2\sum_{l=r+2}^{\nu}u_l \pa_{\zb_l}
+i\sum_{k=\rho +1}^{n-\nu} v_k \pa_{t_k}\in \V^E.
\end{equation}

The basis forms in (\ref{basis}) can be smoothly extended to sections of
$T'M$ near $p$.  Their pullbacks under $b$ are certainly smooth, so
the result follows if we know that these pullbacks remain linearly
independent at $\pt$ and have as common kernel the subspace $\V_{\pt}B$
defined above.  This is straightforward to check using (\ref{basis}),
(\ref{b}), and (\ref{Lform}).
\end{proof}

Denote by $V\subset T\Sigma $ the vertical bundle for 
$b|_\Sigma :\Sigma \rightarrow E$, so that 
$V = \ker b_*|_{T \Sigma}$.   
According to Theorem \ref{involB}, we have 
$\C V\subset \V B|_{\Sigma}$.  If $\pt =(p,[L])\in \Sigma$, then it is
straightforward to calculate from (\ref{b}) that 
$$
b_*(T_{\pt} B) = \R \Im L \oplus T_pE;
$$
this is also easily seen from the geometric interpretation of the blow-up.  
Since $L\in \V_p^E$, it follows that 
\begin{equation}\label{rangeb}
b_*(\C T_{\pt} B) = \C L \oplus \C T_pE.
\end{equation}

If $u$ is a $C^1$ solution of $\V M$, then $u\circ b$ is a solution of 
$\V B$.  (This is also seen to make sense and hold for distribution
solutions, using the smoothness transverse to $E$.)  If $\V M$ is locally
integrable and $\{Z_1,\ldots,Z_m\}$ are solutions with linearly independent
differentials at each point of some set, then it follows from 
$b_*(\C T_{\pt} B)\supset \C T_pE$, $\pt \in \Sigma$, 
and the fact that $E$ is strongly
noncharacteristic, that the differentials of 
$\{Z_1\circ b,\ldots,Z_m \circ b\}$ 
are likewise linearly independent.  Therefore, a hypo-analytic
structure for $(M,\V M)$ naturally lifts to a hypo-analytic structure for
$(B,\V B)$.

If $\pt = (p,[L])\in \Sigma$, we always have 
$\Re L\in T_pE$.  However, it may or may not happen that 
$\Re L \in \Re \V_pE$.  We shall say that $\pt$ is of the first type if 
$\Re L \notin \Re \V_pE$, and that $\pt$ is of the second type if 
$\Re L \in \Re \V_pE$.  In terms of the coordinates introduced in the
proof of Theorem \ref{involB}, lines in $\{(y'',t'')\}$ correspond to
points of the first or second type according to whether $y''$ is not or is
equal to 0.  If $\V M$ is CR, then all points of $\Sigma$ are of the first
type. 

Consider next the characteristic set 
$T^{\circ}_{\pt} B$ at $\pt = (p,[L])\in \Sigma$.  We certainly have 
$T^{\circ}_{\pt} B\subset (V_{\pt})^{\perp}$.  According to (\ref{rangeb}),
any vector $X\in \C T_pE$ can be lifted to a vector 
${\tilde X}\in \C T_{\pt} B$ satisfying $b_*{\tilde X}=X$.  If 
${\tilde \s}\in T^{\circ}_{\pt} B$, then ${\tilde \s}({\tilde X})$ is
independent of the lift and vanishes if $X\in \V _p E$, so we obtain a map  
$T^{\circ}_{\pt} B \rightarrow T^{\circ}_p E$, where, as in the previous
sections, $T^{\circ} E$ denotes the characteristic set for the induced
involutive structure on $E$.  By (\ref{rangeb}), this map is injective.  
We claim that its range is $T^{\circ}_pE\cap (\Re L)^{\perp}$.  Certainly
the range is contained in this set.  Now
any $\s\in T^{\circ}_pE$ has a unique extension to $T_pE\oplus \R\Im L$
annihilating $\Im L$.  This extension
annihilates $L$ if $\s\in (\Re L)^{\perp}$.  According to (\ref{rangeb}), 
this extension may be pulled back under $b^*$ 
to give the required element of $T^{\circ}_{\pt}B$, which with some license
we denote by $b^*\s$.  We therefore have
proved the following:
\begin{proposition}\label{charB}
The pullback $b^*$ induces as described above an isomorphism 
$$
T^{\circ}_{\pt}B \cong T^{\circ}_pE\cap (\Re L)^{\perp}.
$$
\end{proposition}
\noindent
In particular, if $\pt$ is of the first type, then $T^{\circ}_{\pt}B$ is
isomorphic to a codimension 1 subspace of  $T^{\circ}_pE$, while if 
$\pt$ is of the second type, then 
$T^{\circ}_{\pt}B \cong T^{\circ}_pE$. 

Since $\C L \oplus \V_pE\subset \V_pM$, immediately from the definition we
see that $b^*$ defines an injection 
$T^{\circ}_pM\hookrightarrow T^{\circ}_{\pt}B$.  Under the isomorphism of 
Proposition \ref{charB}, this corresponds to the subspace
$\iota_E^*T^{\circ}_pM\subset T^{\circ}_pE\cap (\Re L)^{\perp}$ in the
notation of the previous sections.

One can also easily describe invariantly much of the Levi form of $B$ 
at points of $\Sigma$.  To do so, recall that if $\bV$ is a real 
vector space and $\ell \in \bP(\bV)$ 
is a line in $\bV$, then there is a canonical isomorphism 
$T_{\ell}\bP(\bV)\cong \Hom (\ell,\bV/\ell).$  For $v\in T_{\ell}\bP(\bV)$,   
the corresponding homomorphism 
$\xi _{v}\in \Hom (\ell,\bV/\ell)$ is defined by 
$\xi_v (L) = {\tilde v} +\ell$, where $0\neq L\in \ell$ and 
${\tilde v}\in T_L \bV\cong \bV$ is any tangent vector with 
$\pi_*({\tilde v})=v$, and $\pi:\bV \setminus \{0\}\rightarrow \bP(\bV)$  
denotes the defining projection.  Given $\pt = (p,[L])\in \Sigma$, we apply
this taking $\bV = \V_p^E/\V_p E$ to obtain for
$v\in V_{\pt}$ a homomorphism 
$\xi_v:[L]\rightarrow \V_p^E/(\R L\oplus \V _pE)$.  

Recall from Proposition \ref{charB} that every element of 
$T_{\pt}^{\circ}B$ is of the form $b^*\s$ for some 
$\s\in T^{\circ}_pE\cap (\Re L)^{\perp}$, and from Theorem \ref{involB} 
that $(b_*)^{-1}(\V _p E)$ is a codimension 1 subspace of 
$\V_{\pt}B$.  Direct calculation in
local coordinates gives the following:
\begin{proposition}
Let $\pt = (p,[L])\in \Sigma$ and let 
$\s\in T^{\circ}_pE\cap (\Re L)^{\perp}$.
Then
\begin{enumerate}
\item
$\Le_{b^*{\s}}(L_1,L_2) = \Le_{\s}(b_*L_1,b_*L_2),\,\,$
if either 
\begin{enumerate}
\item $L_1, L_2\in (b_*)^{-1}(\V _p E)$, or
\item $\s \in T^{\circ}_p M$ and $L_1, L_2 \in \V_{\pt}B$.
\end{enumerate}
\item
$\Le_{b^* \s}(L_1,v) = i\s (\Re \xi _v(b_*L_1)),\,\,$ for \\
$L_1\in (b_*)^{-1}(\R L \oplus \V_pE)$ and $v\in V_{\pt}$.
\end{enumerate}
\end{proposition}
\noindent
In 1.(a), the Levi form on the right hand side of 1.\ is that of $E$, 
while in 1.(b) it is that of $M$.  Observe that the right hand side of 1.\
is not defined for general $\s\in T^{\circ}_pE\cap (\Re L)^{\perp}$ and 
$L_1, L_2\in \V_{\pt}B$.  On the right hand side of 2., we 
view $b_*L_1 \in [L]$ as usual as determined up to addition of an
element of $\V_pE$.  Also, $\s$ denotes the induced linear functional on 
$\Re (\V_p^E/(\R L\oplus \V _pE)) \subset T_pE/(\R \Re L + \Re \V_pE)$.
  
{From} either 1.\ or 2., one concludes that 
$\Le_{\sit}(L_1,v)=0$ if $\sit\in T^{\circ}_{\pt}B$, $v\in \C V_{\pt}$,  
and $L_1\in (b_*)^{-1}(\V_p E)$.
In particular, $\Le_{\sit}(v_1,v_2)=0$ for $v_1,v_2\in \C V_{\pt}$.  
(This latter is of course immediate from the fact that $V$ is an integrable 
subbundle of $T\Sigma$.)  Also, it follows from 1.(b) that for
$\s\in T^{\circ}_pM$, one has
$\Le_{b^*\s}(L_1,v)=0$ for $v\in \C V_{\pt}$ and for all 
$L_1\in \V_{\pt}B$.  However, if 
$\sit \in T^{\circ}_{\pt}B\setminus b^*(T^{\circ}_pM)$, then from 2.\ one 
sees that $\Le_{\sit}(L_1,v)\neq 0$ for some $L_1\in \V_{\pt}B$ and
$v\in V_{\pt}$.  Thus for any 
$\sit \in T^{\circ}_{\pt}B\setminus b^*(T^{\circ}_pM)$, it is the case that
$\Le_{\sit}$ is an indefinite Hermitian form on $\V_{\pt}B$.

The blow-up can be used to prove regularity theorems like those of the
previous section.  The geometry is particularly nice when one has solutions 
on two opposite wedges in $M$, with boundary values which agree on the
edge.  The solutions lift to solutions on $B$ defined on opposite sides of
the hypersurface $\Sigma$, and the fact the boundary values agree implies
that one obtains a solution in a full neighborhood of a point
of $\Sigma$.  Therefore the results of \cite{bct} and \cite{chang} can be
applied on $B$, and the conclusions then reinterpreted
on $M$.  Consider for example Corollary \ref{2sidednonchar}.  The relevant 
fact on $B$ is derived above: 
$\Le_{\sit}$ is an indefinite Hermitian form on $\V_{\pt}B$ for any 
$\sit \in T^{\circ}_{\pt}B\setminus b^*(T^{\circ}_pM)$.  According to 
Theorem 6.1 of \cite{bct}, it follows that 
$WF^B_{\pt}(u)\subset b^*(T^{\circ}_pM)$ for any solution $u$ on $B$.  
This can 
be easily reinterpreted on $M$ to give Corollary \ref{2sidednonchar}. 
Similarly, one can prove two-sided versions of Theorems \ref{negev}
and \ref{nondegnull} by passing to $B$.  For Theorem \ref{negev}, one uses  
again Theorem 6.1 of \cite{bct}, but now applied to
$\sit \in b^*(T^{\circ}_pM)$.  For Theorem \ref{nondegnull}, one uses
Chang's theorem on $B$.  
Observe that often the relevant information on $B$
is observed at one higher order than on $M$. 
Results for solutions on one-sided wedges in $M$ 
reduce to the case of a noncharacteristic hypersurface boundary
on $B$, for which Theorem \ref{main} can be applied.  
However, for both one- and two-sided wedges it seems ultimately
more efficient to argue directly on $M$ as in the previous
sections of this paper.

\end{document}